\newif\iffinal
\finalfalse	
\finaltrue	

\documentclass[letterpaper,12pt,reqno]{amsart}
\RequirePackage[utf8]{inputenc}
\usepackage[portrait,margin=2.9cm]{geometry}
\iffinal\else\usepackage[notref,notcite]{showkeys}\fi
\usepackage{mathrsfs,soul}
\usepackage{hyperref}
\usepackage{relsize}
\usepackage[foot]{amsaddr}
\usepackage{amssymb,amsthm,amsfonts,amsbsy,latexsym,dsfont}
\usepackage[textsize=small]{todonotes}       
\usepackage{graphicx}
\usepackage[numeric,initials,nobysame]{amsrefs}
\usepackage{upref,setspace}
\usepackage{xcolor,colortbl}
\usepackage{enumerate}

\usepackage{paralist}

\newenvironment{inparaenumn}{\begin{inparaenum}[\upshape \bfseries (1) ]}{\end{inparaenum}}

\numberwithin{equation}{section}
\numberwithin{figure}{section}
\numberwithin{table}{section}

\newtheorem{thm}{Theorem}[section]
\newtheorem{lem}[thm]{Lemma}
\newtheorem{theorem}[thm]{Theorem}

\newtheorem{prop}[thm]{Proposition}

\newtheorem{conj}[thm]{Conjecture}
\newtheorem{lemma}[thm]{Lemma}

\theoremstyle{definition}
\newtheorem{rem}{Remark}

\newcommand{\ind}{\mathds{1}}
\newcommand{\eps}{\varepsilon}

\newcommand{\set}[1]{\left\{#1\right\}}

\newcommand{\probc}{\stackrel{\mathrm{P}}{\longrightarrow}}

\newcommand{\weakc}{\stackrel{\mathrm{d}}{\longrightarrow}}
\newcommand{\convas}{\stackrel{\mathrm{a.s.}}{\longrightarrow}}

\def\qed{ \hfill $\blacksquare$}

\newcommand{\cA}{\mathcal{A}}\newcommand{\cB}{\mathcal{B}}\newcommand{\cC}{\mathcal{C}}
\newcommand{\cF}{\mathcal{F}}
\newcommand{\cG}{\mathcal{G}}

\newcommand{\cM}{\mathcal{M}}
\newcommand{\cR}{\mathcal{R}}
\newcommand{\cS}{\mathcal{S}}
\newcommand{\cW}{\mathcal{W}}

\newcommand{\vC}{\mathbf{C}}

\newcommand{\vM}{\mathbf{M}}

\newcommand{\vx}{\mathbf{x}}

\newcommand{\mvc}{\boldsymbol{c}}

\newcommand{\mvu}{\boldsymbol{u}}
\newcommand{\mvx}{\boldsymbol{x}}\newcommand{\mvy}{\boldsymbol{y}}
\newcommand{\mvz}{\boldsymbol{z}}

\newcommand{\mvzeta}{\boldsymbol{\zeta}}
\newcommand{\mvxi}{\boldsymbol{\xi}}

\newcommand{\fC}{\mathfrak{C}}

\newcommand{\bD}{\mathbb{D}}\newcommand{\bE}{\mathbb{E}}

\newcommand{\bN}{\mathbb{N}}
\newcommand{\bP}{\mathbb{P}}\newcommand{\bQ}{\mathbb{Q}}\newcommand{\bR}{\mathbb{R}}

\DeclareMathOperator{\E}{\mathbb{E}}
\DeclareMathOperator{\pr}{\mathbb{P}}

\DeclareMathOperator{\var}{Var}

\DeclareMathOperator{\er}{er}

\newcommand{\sss}{\scriptscriptstyle}
\newcommand{\erdos}{Erd\H{o}s-R\'enyi }

\newcommand{\ldown}{l^2_{\downarrow}}
\newcommand{\lpdown}{l^p_{\downarrow}}

\newcommand{\ER}{{\text{ER}}} 
\newcommand{\convd}{\stackrel{d}{\longrightarrow}}

\newcommand{\op}{o_{\sss \mathrm{P}}}

\usepackage{fourier}

\definecolor{aqua}{rgb}{0.0, 1.0, 1.0}
\definecolor{webbrown}{rgb}{.6,0,0}
\definecolor{pinegreen}{rgb}{0.0, 0.47, 0.44}
\definecolor{ultramarineblue}{rgb}{0.25, 0.4, 0.96}
\definecolor{jrnl}{rgb}{0.0, 0.5, 0.0}
\definecolor{lincolngreen}{rgb}{0.11, 0.35, 0.02}
\definecolor{green(html/cssgreen)}{rgb}{0.0, 0.5, 0.0}
\definecolor{airforceblue}{rgb}{0.36, 0.54, 0.66}
\definecolor{azure}{rgb}{0.0, 0.5, 1.0}
\definecolor{bleudefrance}{rgb}{0.19, 0.55, 0.91}
\definecolor{cobalt}{rgb}{0.0, 0.28, 0.67}

\hypersetup{colorlinks=true, linktocpage=true, pdfstartpage=1, pdfstartview=FitV,
	breaklinks=true, pdfpagemode=UseNone, pageanchor=true, pdfpagemode=UseOutlines,
	plainpages=false, bookmarksnumbered, bookmarksopen=true, bookmarksopenlevel=1,
	hypertexnames=true, pdfhighlight=/O,
	urlcolor=webbrown, linkcolor=cobalt, citecolor=jrnl
}

\newcommand{\ch}[1]{\textcolor{black}{{#1}}}
\newcommand{\chch}[1]{\textcolor{black}{{#1}}}

\usepackage{enumitem}
\setlist[itemize]{leftmargin=*}

\begin{document}
\title[Leader problem in random graphs]{A probabilistic approach to the leader problem in random graphs}

\date{}

\author[Addario-Berry]{Louigi Addario-Berry$^1$}
\address{\hskip-15pt $^1$Department of Mathematics and Statistics, McGill University, Montreal, Canada}
\author[Bhamidi]{Shankar Bhamidi$^2$}
\address{\hskip-15pt $^2$Department of Statistics and Operations Research, University of North Carolina, Chapel Hill, USA}
\author[Sen]{Sanchayan Sen$^3$}
\address{\hskip-15pt $^3$Department of Mathematics, Indian Institute of Science, Bangalore, India}
\email{louigi.addario@mcgill.ca, bhamidi@email.unc.edu, sanchayan.sen1@gmail.com}

\maketitle

\begin{abstract}
We study the fixation time of the identity of the leader, i.e., the most massive component, in the general setting of Aldous's multiplicative coalescent \cite{aldous-crit, aldous-limic}, which in an asymptotic sense describes the evolution of the component sizes of a wide array of near-critical coalescent processes, including the classical \erdos process.

We show tightness of the fixation time in the ``Brownian'' regime, explicitly determining the median value of the fixation time to within an optimal $O(1)$ window.
This generalizes {\L}uczak's result \cite{luczak1990component} for the \erdos random graph using completely different techniques.

\chch{
In the heavy-tailed case, in which the limit of the component sizes can be encoded using a thinned pure-jump L\'{e}vy process, we prove that only one-sided tightness holds.
This shows a genuine difference in the possible behavior in the two regimes.
}

The solution to the leader problem in the setting of the \erdos random graph played an important role in the study of the scaling limit of the minimal spanning tree on the complete graph \cite{AddBroGolMie13}. We believe that analogous results, such as those proved herein, will be useful in establishing universality of the intrinsic geometry of the minimal spanning tree across a large class of models.
\end{abstract}

\section{Introduction and main results}\label{sec:intro}
\ch{The foundational work of \erdos \cite{erdos-renyi-1,erdos-renyi-2} motivated an enormous amount of work on the study of dynamically evolving random graph models
in the ensuing decades.} We now briefly describe one of the motivating questions of this paper and then \ch{discuss} renewed interest on this problem over the last few years.

One of the main models studied in the original work of \erdos in \cite{erdos-renyi-1} is the following ``random graph process'' $\big(\ER(n,M),~M\geq 0\big)$ on $[n] := \set{1,2,\ldots, n}$. 
Set $\ER(n,0)$ to be the empty graph. 
For $ M \geq 1$, $\ER(n,M)$ is obtained from $\ER(n,M-1)$ by choosing one of the ${n\choose 2} - M+1$ edges not present in $\ER(n,M)$ uniformly at random and placing this in the system. 
Write $\cC_{\sss(k)}(M)$ (respectively $|\cC_{\sss(k)}(M)|$) for the $k$-th largest component in $\ER(n,M)$ (respectively the size of this component), breaking ties arbitrarily. 
Here we have suppressed dependence on $n$ to simplify notation.
In \cite{erdos-renyi-2}, the following ``double jump'' was identified where it was shown that for $M \ll n/2$\ch{,} $|\cC_{\sss(1)}(M)| = O_P(\log{n})$, if $M = n/2$ \chch{(which corresponds to the so-called {\bf critical regime})}, then $|\cC_{\sss(1)}(M)| = \Theta_P(n^{2/3})$, whilst if $M = c n/2$  with $c> 1$, then $|\cC_{\sss(1)}(M)| \sim f(c) n$ for a deterministic function $f$ satisfying $f(c) > 0$ for $c> 1$. 
This stimulated an enormous amount of work (see \cite{Boll-book,janson2011random,luczak1990component,luczak1994structure} and the references therein) \chch{on understanding the behavior close to the critical regime, and the} dynamic properties of the above construction wherein components merge via the addition of new edges.  This resulted in the following fundamental result of Aldous \cite{aldous-crit}.  Fix $\lambda \in \bR$ and consider the process of normalized component sizes close to the critical value arranged in decreasing order:
\[\bar{\vC}_n(\lambda) =\Big(n^{-2/3}\big|\cC_{\sss(k)}\left(n/2 + \lambda n^{2/3}\right)\big|,~ k\geq 1\Big)\, \chch{.}  \]
For any $p\geq 1$, consider the metric space
\begin{equation}
\label{eqn:ldown-def}
\lpdown:= \bigg\{\vx= (x_i:i\geq 1): x_1\geq  x_2 \geq \ldots \geq 0, \sum_{i=1}^\infty x_i^p < \infty\bigg\}
\end{equation}
equipped with the natural metric inherited from \chch{$l^p$}.
\begin{theorem}[\cite{aldous-crit}]\label{thm:aldous-erdos-ren}
	View the process $\set{\bar{\vC}_n(\lambda): -\infty < \lambda < \infty}$ as a Markov process on $\ldown$. Then as $n\to\infty$, the finite dimensional distributions of the above process converges weakly to that of a Markov process on $\ldown$
	which is referred to as the standard multiplicative coalescent.
\end{theorem}
We will describe this result (as well as the entrance boundary of the Markov process) in more detail in Section \ref{sec:bfs-const}; much more extensive discussions of this process and the relationship to the evolution of the \erdos random graph can be found in \cite{aldous-crit}. We are now in a position to state the main problem motivating this paper.

\noindent {\bf Leader problem:} Erd\H{o}s suggested that one should view the original random graph process $\big(\ER(n,M),~M\geq 0\big)$ as a ``race of components.'' 
One fascinating aspect of this view was studied in \cite{luczak1990component}. First we need some definitions. \chch{For a graph $G$, we call any connected component of $G$ with the maximum number of vertices a leader. 
Now consider the \erdos process $\big(\ER(n,M),~M\geq 0\big)$. 
For $0\leq M_1<M_2$, we say that a {\it change of leader does not occur} in the interval $[M_1, M_2]$ if there exists a leader $\cC$ in $\ER(n, M_1)$ such that for all $M\in [M_1, M_2]$, the component in $\ER(n, M)$ containing $\cC$ is a leader in $\ER(n, M)$.
}
Define
\begin{align}\label{eqn:erdos}
	L^{\sss(n)}_{\er}:=\min\big\{s\geq 0\ : \text{ a change of leader does not occur in the process }\big(\ER(n, M),~M\geq s\big)\big\}.
\end{align}
Then {\L}uczak in \cite[Theorem 7]{luczak1990component} showed that
\begin{equation}
\label{eqn:luczak-tight}
\text{ the sequence of random variables } \set{n^{-2/3}\big(L^{\sss(n)}_{\er}-n/2\big)}_{n\geq 1} \text{ is tight.}
\end{equation}

\noindent {\bf Aim of this paper:} The original proof in \cite{luczak1990component} is highly intricate using careful and refined combinatorial analysis of the number of components of various complexities including trees of various sizes, coupled with a ``symmetry rule'' relating properties of the process below and above the critical threshold. These estimates are combined with a ``scanning method'' to prove \eqref{eqn:luczak-tight}. This paper is motivated by the following two threads:

\begin{inparaenumn}
    \noindent\item 
	In the last few years, a host of random graph models have been shown to belong to the \erdos or more precisely, the multiplicative coalescent universality class \cite{van2013critical,bhamidi2010scaling,SBVHVJL12,bhamidi2014augmented, SBSSXW14,bhamidi-broutin-sen-wang, bhamidi-sen, joseph2014component,riordan2012phase,dhara-hofstad-leeuwaarden-sen}. This includes the configuration model \cite{bollobas1980probabilistic,molloy1995critical}, a large sub-class of the inhomogeneous random graph models as formulated in \cite{BJR07}, and the so-called bounded size rules \cite{spencer2007birth}. It is hard to generalize {\L}uczak's result in \eqref{eqn:erdos} \chch{to these random graph models} via the beautiful counting arguments in \cite{luczak1990component}. Thus the aim of this paper is to give a short probabilistic proof of the above result that is robust and applies to the general setting of the multiplicative coalescent. The classical \erdos case can be recovered from our result; \ch{moreover} the techniques in this paper apply to a number of other entrance boundary conditions for the multiplicative coalescent that have arisen in the study of \emph{heavy-tailed} critical random graphs \cite{SBVHVJL12,joseph2014component,dhara-hofstad-leeuwaarden-sen}.

	\noindent\item 
	Coupled with renewed interest in the critical regime, the last few years have also witnessed an explosion in the study of various models of information propagation on network models. In this context one major model that has been explored in great detail is the minimal spanning tree ({\bf MST}) problem \cite{BraBulCohHavSta03,braunstein2007optimal}. Here one typically starts with a network model in the supercritical regime (having a giant component). Each edge is assumed to have a random positive length sampled in an \emph{i.i.d.} fashion across edges from a continuous distribution on $(0,\infty)$. The aim then is to understand the (metric) structure of the MST on the giant component, e.g., the typical distance between points on the MST. \ch{To} date the only rigorous results in this context are those in \cite{AddBroGolMie13, addarioberry-sen}.
	In \cite{AddBroGolMie13}, the following was shown:
		      Consider the MST $\cM_n$ on the complete graph on $[n]$ and view this as a tree with edge length one. Rescale each edge of $\cM_n$ by $n^{-1/3}$ and view this a compact metric space $\bar{\cM}_n$. Then as $n\to\infty$, $\bar{\cM}_n$ converges in the Gromov-Hausdorff sense to a limiting random compact metric space $\cM_{\sss (\infty)}$.		
		      
	A key ingredient in \ch{the proof of the main result in \cite{AddBroGolMie13} is the following theorem proved in} \cite{BBG-12}: Consider the \erdos random graph in the critical regime as in the setting of Theorem \ref{thm:aldous-erdos-ren}.
	Fix $\lambda\in \bR$ and $k\geq 1$. Consider $\cC_{\sss(k)}(n/2+\lambda n^{2/3})$ as a metric space where each edge has length one.  Write $\bar{\cC}_{\sss(k)}(\lambda)$ for the resulting metric space where each edge is rescaled by $n^{-1/3}$. Then there exist limiting random compact metric spaces $\vM_{\sss(k)}$ such that $\bar{\cC}_{\sss(k)}(\lambda)$  converge in the Gromov-Hausdorff sense to $\vM_{\sss(k)}$.

	The $n^{-1/3}$ scaling \ch{in the results of both \cite{BBG-12} and \cite{AddBroGolMie13}} is not a coincidence. A key step in the proof in \cite{AddBroGolMie13} is showing that the MST $\bar{\cM}_n$ is close (in a strong sense) to 
	the restriction of $\bar{\cM}_n$ to the maximal component in the critical regime $\bar{\cC}_{\sss(1)}(\lambda)$ ``for a large $\lambda$.'' 
	The proof of this statement relies on the leader result of {\L}uczak \cite{luczak1990component} implying that for large $\lambda$, the identity of the maximal component does {\bf not} change. 
	
	It is believed that under some very general assumptions on the underlying discrete structure, the scaling limit of the minimal spanning tree exists and is universal up to some constants.	
	The problem of establishing universality of the scaling limit of the minimal spanning tree is currently open.
	One way to approach the problem of universality that has proven to be useful \cite{bhamidi-broutin-sen-wang, bhamidi-hofstad-sen, SBSSXW14, bhamidi-dhara-hofstad-sen} would be to first obtain the scaling limit of the MST constructed on an inhomogeneous random graph that is closely related to the multiplicative coalescent, and then extending it to more general random graph models (including those with heavy tailed degrees) by suitable coupling methods.
    Carrying out this program requires the extension of the leader result to the setting of the multiplicative coalescent, which is accomplished in this paper.
\end{inparaenumn}

\subsection{The multiplicative coalescent}\label{sec:mult-coal}
The {\em multiplicative coalescent} is a Markov process with state space $\ldown$.
For a state $\mvx=(x_i, i\geq 1)\in\ldown$, $x_i$ represents the weight or mass of the $i$-th largest cluster.
In the evolution of the multiplicative coalescent, each pair $\{i,j\}$ of clusters merges at rate $x_i x_j$; such a merger results in the state obtained by removing $x_i$ and $x_j$ from $\mvx$ and adding a new entry $x_i+x_j$ inserted at the location corresponding to its rank.
We will write
\[(\mvx(t),t \ge 0)= \big((x_1(t), x_2(t),\ldots),\, t\geq 0\big)\]
for the multiplicative coalescent process started from $\mvx(0) = \mvx$; thus $x_1(t)\geq x_2(t)\geq\ldots\geq 0$ for all $t\geq 0$.

The following is a specific construction of the multiplicative coalescent with initial weights $\mvx$.
Define a random graph process $(G(\mvx,t),t \ge 0)$ as follows. 
For $t\geq 0$, the vertices of $G(\mvx, t)$ are the positive integers $\bN$.
$G(\mvx, 0)$ is the empty graph on $\bN$, and for each $i,j \in \bN$, edges between $i$ and $j$ arrive according to a Poisson point process with rate $x_ix_j$.
In this process, two distinct connected components with masses $a$ and $b$ merge at rate $ab$, where mass of a component $\cC$ is $\sum_{i\in\cC}x_i$. 
Thus, we may couple the processes $G(\mvx,\cdot)$ and $\mvx(\cdot)$ so that for all $t \ge 0$, $\mvx(t)=\big(x_i(t),i \ge 1\big)$ is the ordered sequence of masses of connected components of $G(\mvx,t)$.
We work with this construction throughout.

\chch{
For $t\geq 0$, any component in $G(\mvx, t)$ with mass $x_1(t)$ is called a leader in $G(\mvx, t)$.
For $0\leq t_1<t_2$ with $t_2\in\bR_{>0}\cup\{\infty\}$, we say that a {\it change of leader does not occur} in $\big(G(\mvx, t)\, , \ t\in [t_1, t_2)\big)$ if there exists a leader $\cC$ in $G(\mvx, t_1)$ such that for all $t\in(t_1, t_2)$, the component in $G(\mvx, t)$ that contains $\cC$ is a leader in $G(\mvx, t)$.
For $k\geq 1$, we say that the leader changes at least $k$ times in $[t_1, t_2)$ if there exist $t_1=u_1<u_2<\ldots<u_{k+1}=t_2$ such that a change of leader occurs in $[u_j, u_{j+1})$ for $j=1,\ldots, k$.
We will simply write ``a change of leader does (not) occur in $G(\mvx, \cdot)$'' to mean ``a change of leader does (not) occur in $\big(G(\mvx, t),\ t\geq 0\big)$" etc.}
Our results address the behavior of the random times:
\begin{align}\label{eqn:def-leader}
L(\mvx) = \inf\, \big\{s \ge 0: \mbox{ a change of leader does not occur in }\big(G(\mvx, t), ~t \ge s\big)\big\}\, .
\end{align}

For any $\mvz=(z_1,z_2,\ldots)\in[0,\infty)^{\bN}$ and $r\geq 1$, we will write
\begin{align}\label{eqn;def-sigma-r}
\sigma_r(\mvz):=\sum_{i\geq 1}z_i^r \in [0,\infty].
\end{align} 
In this paper, we will work with a sequence $\big(\mvx^{(n)},\, n \ge 1\big)$ of starting configurations, where $\mvx^{(n)}=(x_i^{(n)}, i\geq 1)\in\ldown$ with $x_i^{(n)}>0$ for only finitely many integers $i$.
We assume throughout that there exists $\mvc = (c_i,i \ge 1)$ such that for all $i \ge 1$,
\begin{equation}\label{eq:main_conv}
\frac{x^{(n)}_i}{\sigma_2(\mvx^{(n)})} \longrightarrow c_i \ \ \text{ as }\ \ n\to\infty\, .
\end{equation}
We consider two conditions that correspond to two different regimes. \\
\noindent{\bf Condition I (Brownian limit).} The sequence $\mvc$ has $c_i = 0$ for all $i \ge 1$, and
\[\frac{\sigma_3(\mvx^{\sss(n)})}{\sigma_2(\mvx^{\sss(n)})^3}\to 1\ \mbox{ as }\ n \to \infty \,.
\]
\noindent{\bf Condition II (Pure jump limit).}
The sequence $\mvc$ is such that
\begin{align}\label{eqn:c}
Ki^{-\alpha}\le c_i\le K'i^{-\alpha}\ \text{ for some }\ \alpha \in (1/3,1/2),
\end{align}
where $0<K\le K'$ are universal constants, and
\begin{align}\label{eqn:555}
\frac{\sigma_3(\mvx^{\sss(n)})}{\sigma_2(\mvx^{\sss(n)})^3}\to \sigma_3(\mvc)\ \mbox{ as }\ n \to \infty \,.
\end{align}

It is easy to show that under both Condition I and Condition II, 
\begin{align}\label{eqn:sigma-2-to-zero}
\sigma_2(\mvx^{\sss(n)})\to 0\,\chch{ ,\ \ \ \text{as}\ \ \ n\to\infty\, .}
\end{align}
Indeed, since $\sigma_3(\mvx^{(n)}) \le x^{(n)}_1 \sigma_2(\mvx^{(n)})$, under Condition I we have
\[
\sigma_2(\mvx^{(n)})^2 = \big(1+o(1)\big) \frac{\sigma_3(\mvx^{(n)})}{\sigma_2(\mvx^{(n)})}
\le (1+o(1)) x^{(n)}_1 = o(1)\cdot \sigma_2(\mvx^{(n)})\, ,
\]
which implies \eqref{eqn:sigma-2-to-zero}.
Under Condition II, for any $M \in \mathbb{N}$, by (\ref{eq:main_conv}) we have
	\[
	\sigma_2(\mvx^{(n)})^2  =
	\big(1+o(1)\big) \frac{\sum_{i \le M} (x_i^{(n)})^2}{\sum_{i \le M} c_i^2} \le (1+o(1)) \frac{\sigma_2(\mvx^{(n)})}{\sum_{i \le M} c_i^2}\, .
	\]
	By assumption, $\sigma_2(\mvc)=\infty$, so the above \ch{bound again yields \eqref{eqn:sigma-2-to-zero}. 
	The fact that  $\sigma_2(\mvx^{\sss(n)})\to 0$} will be used repeatedly throughout the paper.

\begin{rem}
In terms of random graphs, Condition I corresponds to the mean-field percolation regime and Condition II corresponds to the heavy-tailed regime.
More precisely, the component sizes of barely-subcritical random graphs whose degree distribution obeys a power law with exponent $\tau=(\alpha+1)/\alpha \in (3,4)$ satisfies \eqref{eq:main_conv} with $\mvc$ as in Condition II; see \cite{van2013critical,SBVHVJL12,joseph2014component,dhara-hofstad-leeuwaarden-sen-2}.
In that setting, in a graph with $n$ nodes, the $i$-th largest degree will typically be of order $(n/i)^{\alpha}$. 
Limiting vectors $\mvc$ of the form \eqref{eqn:c} arise in the description of the scaling limits of such random graphs at criticality.
\end{rem}

\section{Main results}\label{sec:main-results}
We state our main results in this section.
\begin{thm}\label{thm:main_i}
Under Condition I, the collection of random variables
\[
\left\{L(\mvx^{(n)}) - \frac{1}{\sigma_2(\mvx^{(n)})}\, :\, n \ge 1\right\}
\]
is tight.
\end{thm}
This result says that the leader changes for the last time around time $t=1/\sigma_2(\mvx^{(n)})$ with $O(1)$ fluctuations. 
\ch{One can recover \eqref{eqn:luczak-tight} from Theorem \ref{thm:main_i}.
We will elaborate on this in Remark \ref{rem:1}.}

Under Condition II, we show that the leader does not change much after time $1/\sigma_2(\mvx^{(n)})$.
\begin{thm}\label{thm:main_ii}
Under Condition II, the collection of random variables
\[
\left\{\left(L(\mvx^{(n)}) - \frac{1}{\sigma_2(\mvx^{(n)})}\right)^{\mathlarger{+}}\, :\,  n \ge 1\right\}
\]
is tight.
\end{thm}
\ch{There is a genuine difference in the possible behavior in the two regimes as shown in the next theorem.
\begin{thm}\label{thm:iv}
	Under Condition II, whenever $c_1>c_2$,
	\[
	\lim_{\lambda\to\infty}\liminf_{n \to \infty}\
	\pr\bigg(L(\mvx^{(n)})\leq -\lambda+1/\sigma_2(\mvx^{(n)})\bigg)>0\, . 
	\]
\end{thm}
\begin{rem}
It would be interesting to see if the following stronger result holds:
Under Condition II, $\liminf_{n\to\infty} \pr\big(L(\mvx^{(n)})=0\big) > 0$ if $c_1 > c_2$.
If this is not true, then the natural question would be the following:
What is a threshold function $a_n=a_n(\mvx^{(n)})\to\infty$ such that under Condition II, if $c_1 > c_2$, then
\begin{align*}
&\liminf_{n\to\infty}\ \pr\big(L(\mvx^{(n)})\leq -\lambda_n+1/\sigma_2(\mvx^{(n)})\big) > 0\, ,
\ \ \text{ if }\ \ \lambda_n<<a_n\, , \ \ \text{ and}\\
&\lim_{n\to\infty}\ \pr\big(L(\mvx^{(n)})\leq -\lambda_n+1/\sigma_2(\mvx^{(n)})\big) = 0\, ,
\ \ \text{ if }\ \ \lambda_n>>a_n\, .
\end{align*}
\end{rem}}
In contrast,  under Condition I, the following result holds.
\begin{thm}\label{thm:number-of-changes}
Let $N^{\sss(n)}$ denote the number of times the leader changes in $\big[0, \big(\sigma_2(\mvx^{\sss(n)})\big)^{-1}\big]$. Under Condition I, $N^{\sss(n)}\probc\infty$.
\end{thm}
\ch{
It is known \cite{aldous-crit, aldous-limic} that under both Condition I and Condition II, there exists an $\ldown$-valued process $(\mvzeta(\lambda), \lambda\in\bR)$ (called an eternal multiplicative coalescent) such that for every fixed $\lambda\in\bR$, $\mvx^{(n)}(t_{\lambda})\convd \mvzeta(\lambda)$ as $n\to\infty$, where $t_{\lambda}=t_\lambda(\mvx^{(n)})=\lambda+1/\sigma_2(\mvx^{\sss(n)})$; this is discussed in more detail in Section \ref{sec:bfs-const}. (Theorem \ref{thm:aldous-erdos-ren} is a special case of this result.)
The arguments in the proofs of Theorems \ref{thm:main_i}, \ref{thm:main_ii}, \ref{thm:iv}, and \ref{thm:number-of-changes} carry over in an identical way to the continuum yielding analogous results for $\mvzeta(\cdot)$.
We collect these results in the next theorem; we omit the proof.
\begin{theorem}\label{thm:iii}
{\upshape (a)}
If $c_i=0$ for all $i\geq 1$, or if $\mvc$ satisfies \eqref{eqn:c}, then 
\[
\lim_{\lambda\to\infty}\ \pr\bigg(\text{a change of leader occurs in }G\big(\mvzeta(\lambda),\ \cdot\big)\bigg)=0\, .
\]
{\upshape (b)} 
If $c_i=0$ for all $i\geq 1$, then for all $k\geq 1$,
\[
\lim_{\lambda\to \infty} \pr\left(\text{the leader changes at least }k\text{ times in }\bigg(G\big(\mvzeta(-\lambda), \ t\big)\, , \ t\in[0, \lambda]\bigg)\right)=1\, .
\]
{\upshape (c)} 
If $\mvc$ satisfies \eqref{eqn:c}, then 
\[
\lim_{\lambda\to \infty} \pr\bigg(\text{a change of leader does not occur in }G\big(\mvzeta(-\lambda),\ \cdot\big)\bigg)>0\, .
\]
\end{theorem}
Now, Theorems \ref{thm:main_i} and \ref{thm:main_ii} state results about tightness.
A natural question about convergence in distribution arises here.
In this context, we make the following conjecture.
\begin{conj}\label{conj:2}
Under both Condition I and Condition II, for every $\lambda\in\bR$,
\[
\lim_{n\to\infty}\pr\Big(L(\mvx^{(n)})> \big(\sigma_2(\mvx^{(n)})\big)^{-1}+\lambda\Big)
=
\pr\bigg(\text{a change of leader occurs in }G\big(\mvzeta(\lambda), \ \cdot\big)\bigg)\, .
\]
\end{conj}
}

\begin{rem}
In \cite{luczak1990component}, {\L}uczak studied the size of the leader in the setting of the \erdos random graph.
Informally, his results (in particular \cite[Theorem 5]{luczak1990component}) give a law of large numbers for the mass of the largest component at time $1/\sigma_2 + \lambda$, as $\lambda \to \infty$, in the Erd\H{o}s-R\'enyi setting. {\L}uczak's results in particular imply that under Condition I, the mass of the leader at time $1/\sigma_2 + \lambda$ grows linearly in $\lambda$. Our Theorem~\ref{lem:largest-component-big} stated below, provides an analogous result under Condition II. Under Condition II, the growth rate is not $\lambda$ but $\lambda^{(1-2\alpha)/\alpha}$. Related results for the barely supercritical configuration model with heavy-tailed degrees were proved in \cite{hofstad16configuration}.
\end{rem}

\begin{rem}\label{rem:1}
\ch{Let us briefly describe how to prove \eqref{eqn:luczak-tight} using Theorem \ref{thm:main_i}.
Consider the modified \erdos process $\big(\overline\ER(n, M),\ M\geq 0\big)$, where at each step, edges are sampled uniformly and with replacement, i.e., multiple edges are allowed.
A simple computation shows that the number of multiple edges created up to $M\leq n$ is tight.
Thus, it is enough to show that $n^{-2/3}\big(\overline L_{\er}^{\sss(n)}-n/2\big)$ is tight, where $\overline L_{\er}^{\sss(n)}$ denotes the minimum of all $s\geq 0$ such that a change of leader does not occur in $\big(\overline\ER(n, M),\ M\geq s\big)$.}

\ch{
Consider the weight sequence 
$\mvy^{(n)}=(y^{(n)}_i,~i\geq 1)$, where 
\begin{equation}\label{eq:errg-weights}
y^{(n)}_i =
\begin{cases}
n^{-2/3},&\text{ if }\ i \le n \\
0,& \text{ if }\ i > n.
\end{cases}
\end{equation}
Then $\mvy^{\sss(n)}$ satisfies Condition I and $\sigma_2(\mvy^{\sss(n)})=n^{-1/3}$.
Consider the random graph process $\big(G(\mvy^{\sss(n)},t),\ t\geq 0\big)$, and for $k\geq 1$, let $\tau_k$ be the time when the $k$-th edge is added in this process.
Clearly $\big(\overline\ER(n, M),\ M\geq 1\big)$ has the same distribution as the process $\big(G(\mvy^{\sss(n)},\tau_M),\ M\geq 1\big)$ restricted to the vertex set $\{1,\ldots, n\}$.
For $u\in\bR$, define $M_n(u)=[n/2+un^{2/3}]$.
(Note that for any $u$, $M_n(u)>0$ for all large $n$.)
Then,
\begin{align}\label{eqn:43}
\bP\big(\overline L_{\er}^{\sss(n)}\geq M_n(u)\big)
=\bP\big(L(\mvy^{\sss(n)})\geq\tau_{M_n(u)}\big).
\end{align}
Now, edges appear at a rate of $\lambda_n:=n^{-4/3}{n\choose 2}$ in the process $G(\mvy^{\sss(n)}, \cdot)$.
Hence, for any fixed $u\in\bR$,
\begin{align*}
\tau_{M_n(u)}
=\lambda_n^{-1}M_n(u)+O_P\big(\lambda_n^{-1}\sqrt{M_n(u)}\big)
=n^{1/3}+2u+O_P(n^{-1/6})
=\big[\sigma_2(\mvy^{\sss(n)})\big]^{-1}+2u+O_P(n^{-1/6})\, .
\end{align*}
(Here $O_P(a_n)$ represents a sequence $\big(Y_n,\ n\geq 1\big)$ of random variables such that $\big(a_n^{-1} Y_n,\ n\geq 1\big)$ is tight.)
So \eqref{eqn:43} together with Theorem \ref{thm:main_i} shows that $n^{-2/3}\big(\overline L_{\er}^{\sss(n)}-n/2\big)$ is tight.}
\end{rem}



\section{Preliminaries}
\label{sec:prelim}
\subsection{Notation}
\label{sec:not}
Throughout this paper, we make use of the following notation.
We will write $\convd$, $\probc$, and \ch{$\convas$} to denote convergence in distribution, convergence in probability, and almost sure convergence respectively.
For a sequence of random variables $\big(X_n,\ n\geq 1\big)$, we write $X_n=\op(b_n)$ when $|X_n|/b_n\probc 0$ as $n\rightarrow\infty$.
For a non-negative function $n\mapsto g(n)$, we write $f(n)=O(g(n))$ when $|f(n)|/g(n)$ is uniformly bounded, and $f(n)=o(g(n))$ when $\lim_{n\rightarrow \infty} f(n)/g(n)=0$.
Furthermore, we write $f(n)=\Theta(g(n))$ if $f(n)=O(g(n))$ and $g(n)=O(f(n))$.

Throughout this paper, $C, C', K, K'$ will denote positive constants that depend only on $\mvc=(c_i, i\geq 1)$, and their values may change from line to line.
Special constants will be indexed, e.g., $K_{\ref{lem:bound-expectation}}$ etc.
Given two functions $f_1,f_2: [0,\infty)\to [0,\infty)$, we write $f_1\asymp f_2$ on \ch{$A\subseteq [0,\infty)$} if \ch{there exist $0<K'\leq K<\infty$ such that}
\[K' f_2(x)\leq f_1(x)\leq K f_2(x)\ \text{ for all }\ x\in A.\]
Similarly, for two sequences $\{a_n\}_{n\geq 1}$ and $\{b_n\}_{n\geq 1}$, we will write $\{a_n\}\asymp\{b_n\}$ or simply $a_n\asymp b_n$ to mean that \ch{there exist $0<K'\leq K<\infty$ such that}
\[K' b_n\leq a_n\leq K b_n\ \text{ for all }\ n\geq 1\, .\]

\subsection{Tail bounds for sums and suprema}
We will use the following bound on the suprema of collections of observables of empirical \ch{processes}, which is a specialization of  \cite[Theorem 1.1 (b)] {klein-rio05}(see also \cite[Theorem 12.9]{concentrationbook}).
\begin{thm}[\cite{klein-rio05}] \label{thm:empirical_tail}
Let $X_i$, $i \ge 1$, be independent real random variables. Fix $A > 0$ and $n \ge 1$\ch{,} and let $\{f_{i,s}\, :\, 1 \le i \le n,s \in S\}$ be a countable collection of measurable functions from $\bR$ to $[-A,A]$.
Write
\[
Z = \sup\big\{ \sum_{i=1}^n f_{i,s}(X_i)\ :\ s \in S\big\}\, ,
\]
and let
\[
V = 2\, \bE[Z] + \sup_{s \in S} \var\big(\sum_{i \le n} f_{i,s}(X_i)\big)\, .
\]
Then
\[
\bP\bigg(Z \ge \bE Z + x \bigg)
\le
\exp\left(
- \frac{x}{4A} \log\bigg(1 + 2\log\Big(1+\frac{x}{AV}\Big)\bigg)
\right)\, .
\]
\end{thm}
\subsection{Random graph constructions}\label{sec:222}
We work with the construction of the multiplicative coalescent introduced in Section \ref{sec:mult-coal}. 
Recall that in this construction, $\mvx(t)=(x_i(t),i \ge 1)$ is the ordered sequence of masses of the connected components of random graph $G(\mvx,t)$ defined in Section \ref{sec:mult-coal}.
We always assume
\[
\mvx =(x_i,i \ge 1) = (x_i(0),i\ge 1)=\mvx(0) \in \ldown\, .
\]
In this case, $(\mvx(t),t \ge 0)$ is a Markov process with state space $\ldown$. Its generator $\cA$ is given by
\begin{align}\label{eqn:generator}
\cA f(\mvz)=\sum_i\sum_{j>i}z_i z_j\big(f(\mvz^{i,j})-f(\mvz)\big)\, ,
\end{align}
where $\mvz^{i,j} \in \ldown$ is formed from $\mvz$ by removing entries $z_i$ and $z_j$ and inserting a new entry $z_i+z_j$ in a location which preserves membership in $\ldown$.

Let $\cS$ denote the space of all vertex-weighted multigraphs with vertex set $\bN$ that satisfy
\[\sum_{\cC\text{ component of }G}\big(\sum_{i\in \cC}w_i(G)\big)^2<\infty,\]
where $w_i(G)$ denotes the weight of vertex $i$ in $G$.
Then $(G(\mvx,t),t\ge 0)$ is a Markov process with state space $\cS$. Its generator, also denoted $\cA$, is given by
\begin{align}\label{eqn:generator2}
\cA f(G)=\sum_i\sum_{j>i}w_i(G) w_j(G)\big(f(G^{i,j})-f(G)\big)\, ,
\end{align}
where $G^{i,j}$ is formed from $G$ by adding an edge between vertices $i$ and $j$.

For $G\in\cS$ and any connected component $\cC$ of $G$, we write $\cW(\cC) = \sum_{i \in \cC} w_i(G)$ for the mass of $\cC$.
We write $\cC_i(G)$ for the connected component of $G$ having the $i$-th largest mass.
We also write $\cC(G; i)$ for the connected component of $G$ containing vertex $i$.
To simplify notation, we write $\cC_i(\mvx,t)$ and $\cC(\mvx,t; i)$ for the corresponding objects for $G(\mvx,t)$.
Thus, $x_i(t) = \cW(\cC_i(\mvx,t))$ for all $i\geq 1$ and $t\geq 0$, and
$
\sigma_p(\mvx(t)) = \sum_{\cC~\mathrm{component~of}~G(\mvx,t)} \big(\cW(\cC)\big)^p
$
for all $t\ge 0$ and $p \ge 1$.

\subsection{Convergence of multiplicative coalescents}\label{sec:bfs-const}
Define a process $(V(s),s \ge 0)$ as follows.
Under Condition I, the limit sequence $\mvc$ is identically zero. In this case let
\[V(s) = B(s) - s^2/2,\]
where $(B(s), s \ge 0)$ is a standard Brownian motion.
Under Condition II, let $\xi_1,\xi_2,\ldots$ be independent random variables with $\xi_i\sim\text{Exp}(c_i)$,
and let
\begin{align}\label{eqn:2225}
V(s)=
\sum_{i\geq 1}c_i\big(\ind_{\{\xi_i\leq s\}}-c_i s\big)\, .
\end{align}
In both cases, define a process $(W_{\lambda}(s),s \ge 0)$ by setting
\begin{equation}
\label{eqn:wvs-def}
	W_{\lambda}(s)=\lambda s+V(s).
\end{equation}
Let $\overline W_{\lambda}(\cdot)$ denote the process $W_{\lambda}$ reflected at zero, namely,
\begin{equation}
\label{eqn:bar-lam-def}
	\overline W_{\lambda}(s) = W_\lambda(s) - \inf_{0\leq u\leq s} W_\lambda(u), \qquad s\geq 0.
\end{equation}

\ch{Suppose $\mvx\in\ldown$ has finite length,  i.e., has only finitely many nonzero entries.}
A standard tool in the study of the random graph $G(\mvx, t)$ is the ``breadth-first walk" process $\big(B_{\mvx, t}(u), u\geq 0\big)$ associated with a breadth-first exploration of the random graph $G(\mvx, t)$.
\ch{We recall this process from \cite{aldous-limic} briefly; the reader is referred to \cite[Section 2.3]{aldous-limic} for a more detailed description.
Let $U_{i, j}$ $1\leq i, j\leq n$, be independent random variables with $U_{i, j}\sim\text{Exponential}(tx_j )$.
Choose $v(1)$ by size-biased sampling, i.e., vertex $v$ is chosen with probability proportional to $x_v$. 
Define $\{v\, :\, U_{v(1),v}\leq x_{v(1)}\}$ to be the set of children of $v(1)$, and order these children as $v(2), v(3),\ldots$ so that $U_{v(1),v(i)}$ is increasing.
Set $B_{\mvx, t}(0)=0$, and let
\[
B_{\mvx, t}(u)=-u+\sum_v x_v\ind_{\{U_{v(1), v}\leq u\}}\, ,\ \ \ \ 0\leq u\leq x_{v(1)}\, .
\]
Inductively, write $\tau_{i}=\sum_{j\leq i}x_{v(i)}$.
If $v(i)$ is in the same component as $v(1)$, then the set
$
\big\{v\notin\{v(1),\ldots, v(i-1)\}\, :\, v\text{ is a child of one of } \{v(1),\ldots, v(i-1)\}\big\}
$
consists of $v(i),...,v(l(i))$ for some $l(i)\geq i$.
Let the children of $v(i)$ be $\big\{v\notin \{v(1),...,v(l(i))\}\, :\, U_{v(i),v}\leq x_{v(i)}\big\}$, and order them as $v(l(i) + 1), v(l(i)+2),\ldots$ such that $U_{v(i),v}$ is increasing. 
Set
\begin{align}\label{eqn:6666}
B_{\mvx, t}\big(\tau_{i-1}+u\big)
=
B_{\mvx, t}\big(\tau_{i-1}\big)-u+\sum_{v\text{ child of }v(i)} x_v\ind_{\{U_{v(i), v}\leq u\}}\, , \ \ \ \
0\leq u\leq x_{v(i)}\, .
\end{align}
Once the exploration of the component of $v(1)$ is complete, choose the next vertex from the set of the remaining vertices by size-biased sampling, and continue.
In the end, this produces a random forest. 
The partition of the set of vertices into different connected components in this forest has the same distribution as the partition of the set of vertices into different connected components in $G(\mvx, t)$.
Let $\xi_j^n$ denote the `birth-time' of $j$, i.e., the time when vertex $j$ appears in the above exploration.
Thus, if $j$ is the first vertex in its component to appear in the above process, then $\xi_j^n$ is the time when the exploration of that component started.
If $j$ is not the first vertex in its component to appear in the above process and $j$ is the child of $v(i)$, then $\xi_j^n=\tau_{i-1}+U_{v(i), j}$.
}

Now consider $\mvx^{(n)}$, $n\geq 1$, satisfying either Condition I or Condition II.
Now and henceforth, for $\lambda \in \bR$, we write
\[
t_{\lambda}=t_\lambda(\mvx^{(n)})=\lambda+1/\sigma_2(\mvx^{\sss(n)}).
\]
(Note that \eqref{eqn:sigma-2-to-zero} implies that for every $\lambda\in\bR$, $t_{\lambda}(\mvx^{(n)})\geq 0$ for all $n$ sufficiently large.)
Then under both Conditions I and II on the sequence of starting configurations $\big(\mvx^{(n)},\, n \ge 1\big)$, by \cite[Proposition 9]{aldous-limic},
\begin{align}\label{eqn:2928}
\Bigg(\frac{B_{\mvx^{\sss(n)},t_\lambda}(s)}{\sigma_2(\mvx^{(n)})},\, s\geq 0\Bigg)
\weakc
\Big(W_{\lambda}(s),\, s\geq 0\Big).
\end{align}
Let us now look at implication of this convergence.
For an excursion $\gamma$ of $\overline{W}_{\lambda}$, write $|\gamma|$ for the length of $\gamma$.
Let $\gamma_i(\lambda)$, $i\geq 1$, denote the excursions of $\overline W_{\lambda}$ from zero in decreasing order of length, i.e., $|\gamma_1(\lambda)|>|\gamma_2(\lambda)|>\ldots$.
Then under both Condition I and Condition II, by \cite[Proposition 7]{aldous-limic} (see also \cite[Proposition 4]{aldous-crit}), for each fixed $\lambda \in \bR$,
$\mvzeta(\lambda):=(|\gamma_i(\lambda)|,\ i\geq 1)\in\ldown$ almost surely, and further,
\begin{equation}
\label{eqn:weight-tau-con}
	\mvx^{\sss(n)}(t_\lambda)
	\convd
	\mvzeta(\lambda),
\end{equation}
with respect to the topology on $\ldown$.

For any component $\cC$ of $G(\mvx^{(n)}, t_{\lambda})$, let $l_n(\cC)$ (resp. $r_n(\cC)$) denote the time when the exploration of $\cC$ started (resp. concluded) in $B_{\mvx^{\sss(n)},t_\lambda}$.
For $u>0$, let $\fC^{(n)}(\lambda; u)$ denote the component of $G(\mvx^{(n)}, t_{\lambda})$ being explored by $B_{\mvx^{\sss(n)},t_\lambda}$ at time $u$.
For an excursion $\gamma$ of $\overline W_{\lambda}$, write $l(\gamma)$ and $r(\gamma)$ for the start and finish times of $\gamma$.

\ch{
For the rest of Section \ref{sec:bfs-const}, we work under Condition II.
It follows from \cite[Proposition 14 (b)]{aldous-limic} that for every $u>0$, 
$\pr\big(\overline W_{\lambda}(u)=0\big)=0$.
Since $\overline W_{\lambda}(\cdot)$ is continuous at $u$ with probability one for every $u>0$, we conclude that for all $u>0$,
\begin{align}\label{eqn:2227}
\pr\big(\overline W_{\lambda}(u)=0\ \ \text{ or }\ \ \overline W_{\lambda}(u-)=0\big)=0\, .
\end{align}
Consequently, for every $u>0$, 
\begin{align}\label{eqn:2228}
\pr\big(l(\gamma)<u<r(\gamma)\ \text{ for some excursion }\ \gamma\ \text{ of }\ \overline{W}_{\lambda}\big)=1\, .
\end{align}
Denote this excursion by $\gamma(\lambda; u)$.
Now \eqref{eqn:2928} and \cite[Proposition 14~(d)]{aldous-limic} imply (see \cite[Page 27]{aldous-limic}) that the following convergence of point processes holds jointly with \eqref{eqn:2928}:
\begin{align}\label{eqn:pp-convgnce}
\Big\{\big(l_n(\cC), \cW(\cC)\big)\, :\, \cC\text{ component of }G(\mvx^{(n)}, t_{\lambda})\Big\}
\weakc
\Big\{\big(l(\gamma),|\gamma|\big)\, :\, \gamma\text{ excursion of }\overline{W}_{\lambda}\Big\},
\end{align}
the underlying notion of convergence being that of vague convergence of counting measures on $[0, a]\times [b,c]$ for every $a>0$ and $c>b>0$. 
It follows from \eqref{eqn:pp-convgnce} and \eqref{eqn:2228} that for every $u>0$,
\begin{align}\label{eqn:2929-I}
\Big(l_n\big(\fC^{(n)}(\lambda; u)\big),~r_n\big(\fC^{(n)}(\lambda; u)\big)\Big)
\weakc
\Big(l\big(\gamma(\lambda; u)\big),~r\big(\gamma(\lambda; u)\big)\Big)\, 
\end{align}
jointly with the convergence in \eqref{eqn:2928}.
}

\ch{
Recall from \eqref{eqn:2225} and \eqref{eqn:wvs-def} the definitions of the processes $V(\cdot)$ and $W_{\lambda}(\cdot)$ under Condition II. 
Fix $i\geq 1$.
Define
\[
W_{\lambda}^{\sss i-}(s)
=W_{\lambda}(s)-c_i\ind_{\{\xi_i \le s\}}
=\big(\lambda-c_i^2\big)s+\sum_{\substack{j\geq 1\\ j\neq i}}c_j\big(\ind_{\{\xi_j\leq s\}}-c_j s\big)\, ,\ \ \ s\geq 0\, .
\]
Let 
$\overline{W_{\lambda}^{\sss i-}}(s) = W_{\lambda}^{\sss i-}(s) - \inf_{0\leq u\leq s} W_{\lambda}^{\sss i-}(u)$, $s\geq 0$.
Since $\xi_i$ is independent of the process $W_{\lambda}^{\sss i-}(\cdot)$,
the argument leading to \eqref{eqn:2227} shows that 
$
\pr\big(\overline{W_{\lambda}^{\sss i-}}(\xi_i-)=0\big)=0
$.
Using the fact that
$W_{\lambda}^{\sss i-}(s)=W_{\lambda}(s)$ for $s\in[0, \xi_i)$, we see that
$\overline W_{\lambda}(\xi_i-)=\overline{W_{\lambda}^{\sss i-}}(\xi_i-)>0$ with probability one.
Since $\overline W_{\lambda}(\xi_i)>0$, we get
\begin{align}\label{eqn:2229}
\pr\big(\forall j\geq 1\, ,\ \ l(\gamma)<\xi_j<r(\gamma)\ \text{ for some excursion }\ \gamma\ \text{ of }\ \overline{W}_{\lambda}\big)=1\, .
\end{align}
For an excursion $\gamma$ of $\overline W_{\lambda}$ and $j \in \bN$, write $j \in \gamma$ if $l(\gamma)< \xi_j < r(\gamma)$.
}

\ch{
Consider $\eps>0$ such that $c_i>\eps>c_{i+1}$ for some $i\geq 1$.
Then \eqref{eqn:2928}, \eqref{eqn:2929-I}, and \eqref{eqn:2229} imply that for any such $\eps>0$,
\begin{align}\label{eqn:2230}
\sum_{j}\Bigg(\frac{x_j^{(n)}}{\sigma_2(\mvx^{(n)})}\Bigg)^2
\ind\big\{j\in\fC^{(n)}(\lambda; u)\, ;\  x_j^{(n)}>\eps\sigma_2(\mvx^{(n)})\big\}
\weakc
\sum_{j\geq 1}c_j^2\cdot\ind\big\{j\in\gamma(\lambda; u)\, ;\ c_j>\eps\big\}
\end{align}
jointly with the convergence in \eqref{eqn:2929-I}.
(The convergence in \eqref{eqn:2230} can be deduced using the following elementary fact: 
Suppose $f_n\to f$, as $n\to\infty$, in the Skorohod $J_1$ topology on $\bD[0,T]$. 
Then there exists a sequence $\big(\lambda_n,\ n\geq 1\big)$ of strictly increasing, continuous functions on $[0, T]$ satisfying 
(i) $\lambda_n(0)=0$, and $\lambda_n(T)=T$ for all $n\geq 1$, and 
(ii) $\sup_{t\in[0, T]}|\lambda_n(t)-t|+\sup_{t\in[0, T]}|f_n\circ\lambda_n(t)-f(t)|\to 0$ as $n\to\infty$.
In particular, $f_n\circ\lambda_n(t)-f_n\circ\lambda_n(t-)\to f(t)-f(t-)$, as $n\to\infty$, for all $t\in[0, T]$.)
}

\ch{
Let $\xi_j^{n}$ be the birth time of $j$ in the breadth-first walk process as defined below \eqref{eqn:6666}.
As explained in \cite[Page 17]{aldous-limic}, there exist independent random variables $\widetilde\xi_j^n$, $j=1,\ldots, n$, such that $\widetilde\xi_j^n\sim\mathrm{Exponential}(t_{\lambda}x_j^{(n)})$, and $\xi_j^n=\widetilde\xi_j^n$ if $j$ is not the first vertex in its component appearing in the breadth-first walk.
Consequently,
\begin{align}\label{eqn:2232}
\sum_{j\geq k}
\Bigg(\frac{x_j^{(n)}}{\sigma_2(\mvx^{(n)})}\Bigg)^2
\ind\big\{j\in\fC^{(n)}(\lambda; u)\big\}
\leq 
\Bigg(\frac{x_k^{(n)}}{\sigma_2(\mvx^{(n)})}\Bigg)^2
+
\sum_{j\geq k}
\Bigg(\frac{x_j^{(n)}}{\sigma_2(\mvx^{(n)})}\Bigg)^2
\ind\big\{\widetilde\xi_j^n\leq r_n\big(\fC^{(n)}(\lambda; u)\big)\big\}\, .
\end{align}
Now, for any $A>0$,
\begin{align}\label{eqn:2233}
&\lim_{k\to\infty}\limsup_{n \to \infty}\
\bE\Bigg[\sum_{j\geq k}
\Bigg(\frac{x_j^{(n)}}{\sigma_2(\mvx^{(n)})}\Bigg)^2
\ind\big\{\widetilde\xi_j^n\leq A\big\}\Bigg]
\leq
\lim_{k\to\infty}\limsup_{n \to \infty}\
\Bigg[\sum_{j\geq k}
\Bigg(\frac{x_j^{(n)}}{\sigma_2(\mvx^{(n)})}\Bigg)^2 \cdot At_{\lambda}x_j^{(n)}\Bigg]\notag\\
&\hskip40pt=
\lim_{k\to\infty}\limsup_{n \to \infty}\ A\big(1+\lambda\sigma_2(\mvx^{(n)})\big)
\Bigg[\sum_{j\geq k}
\Bigg(\frac{x_j^{(n)}}{\sigma_2(\mvx^{(n)})}\Bigg)^3\Bigg]
=\lim_{k \to \infty}A\sum_{j\geq k}c_j^3=0\, ,
\end{align}
where the penultimate step uses \eqref{eqn:sigma-2-to-zero}, \eqref{eqn:555}, and \eqref{eq:main_conv}, and the final step uses the fact $\mvc\in l^3$.
Since $r_n\big(\fC^{(n)}(\lambda; u)\big)\leq u+ x_1^{\sss (n)}(t_{\lambda})$,
it follows from \eqref{eqn:weight-tau-con} that for every $u>0$, the sequence $\big(r_n\big(\fC^{(n)}(\lambda; u)\big)\, ,\ n\geq 1\big)$ is tight.
Combining this observation with \eqref{eqn:2232} and \eqref{eqn:2233}, we see that for every $\delta>0$,
\begin{align}\label{eqn:2231}
\lim_{k\to\infty}\limsup_{n\to\infty}\ 
\pr\left(
\sum_{j\geq k}
\Bigg(\frac{x_j^{(n)}}{\sigma_2(\mvx^{(n)})}\Bigg)^2
\ind\big\{j\in\fC^{(n)}(\lambda; u)\big\} >\delta
\right)
=0\, .
\end{align}
We can similarly show that 
$\pr\big(\sum_{j\geq k}c_j^2\cdot\ind\big\{j\in\gamma(\lambda; u)\big\}>\delta\big)\to 0$ as  $k\to\infty$, which together with \eqref{eqn:2231} and \eqref{eqn:2230} yields that for every $u>0$ and $\lambda\in\bR$,
\begin{align}\label{eqn:2929-II}
\sum_{j\in\fC^{(n)}(\lambda; u)}\Bigg(\frac{x_j^{(n)}}{\sigma_2(\mvx^{(n)})}\Bigg)^2
\weakc
\sum_{j\geq 1}c_j^2\cdot\ind\big\{j\in\gamma(\lambda; u)\big\}
\end{align}
jointly with the convergence in \eqref{eqn:2929-I}.
}

\section{Proof of Theorem~\ref{thm:main_i}}\label{sec:main_i_proof}
\ch{We break up the claim in Theorem \ref{thm:main_i} into two propositions that we will prove separately.
\begin{prop}\label{prop:4.1}
Under Condition I, 
$\bigg\{\bigg(L(\mvx^{(n)}) - \frac{1}{\sigma_2(\mvx^{(n)})}\bigg)^{\mathlarger{+}}\, :\,  n \ge 1\bigg\}$
is tight.
\end{prop}
\begin{prop}\label{prop:4.2}
Under Condition I, 
$\bigg\{\bigg(L(\mvx^{(n)}) - \frac{1}{\sigma_2(\mvx^{(n)})}\bigg)^{\mathlarger{-}}\, :\,  n \ge 1\bigg\}$
is tight.
\end{prop}
}
The heart of the proof that the family of random variables in Proposition \ref{prop:4.1} is tight is given in the following proposition.
\begin{prop}\label{prop:submartingale}
Fix $\mvy=(y_i,i \ge 1) \in \ldown$. Then for every $\ell\geq 1$,
$\big(f_{\ell}\big(G(\mvy, t)\big), t\geq 0\big)$ is a submartingale, where
\[
f_{\ell}(G)=\frac{\big(\cW(\cC(G; \ell))\big)^2}{\sum_{\cC\text{ component of }G}\cW(\cC)^2},\ \ G\in\cS.
\]
Consequently, if $y_1^2/\sigma_2(\mvy)> 1-\eps$, then
\begin{align}\label{eqn:leader-doesnot-change}
\bP\big(\text{a change of leader does not occur in }G(\mvy, \cdot)\big) > 1-5\eps.
\end{align}
\end{prop}
\noindent{\bf Proof:}
Fix $\mvz\in\ldown$ such that only finitely many coordinates of $\mvz$ are non-zero.
Let $G\in\cS$ be such that $\cW(\cC_i(G))=z_i$ for all $i\geq 1$.
Suppose $\ell\in\cC_k(G)$.
For convenience, we will write $\sigma_2$ for $\sigma_2(\mvz)$.
\ch{We claim that
\begin{equation}\label{eq:submart_toprove}
\cA f_{\ell}(G) \ge z_k\sum_{i\geq 1,i \ne k}\frac{z_i^3}{\sigma_2+2z_iz_k}\, .
\end{equation}
We will prove \eqref{eq:submart_toprove} shortly, but let us first note that since the right-hand side is nonnegative, it follows that the process 
\begin{align}\label{eqn:2234}
\big(f_{\ell}\big(G(\mvy^{[n]}, t)\big), t\geq 0\big)\ \ \text{ is a submartingale, }
\end{align}
where $\mvy^{[n]}=(y_1,\ldots,y_n,0,0,\ldots)$.
To prove the claim for a general $\mvy\in\ldown$, 
note that the graphs $G(\mvy, t)$ and $G(\mvy^{[n]}, t)$ can be coupled in an obvious way so that the latter is a subgraph of the former.
By \cite[Corollary 18 (c)]{aldous-crit}, in this coupling,
\begin{align}\label{eqn:2235}
\sigma_2\big(\mvy^{[n]}(t)\big)\convas \sigma_2\big(\mvy(t)\big)\ \ \text{ as }\ \ n\to\infty\, .
\end{align}
Now, 
$\cW\big(\cC\big(\mvy,t; \ell\big)\big)
=\sum_j y_j\cdot\ind_{\{j\in\cC(\mvy, t; \ell)\}}$.
Since $j\in\cC(\mvy, t; \ell)$ if and only if there is a finite path between $j$ and $\ell$ in $G(\mvy, t)$, $\ind_{\{j\in\cC(\mvy, t; \ell)\}}$ is the almost sure increasing limit of $\ind_{\{j\in\cC(\mvy^{[n]}, t; \ell)\}}$ in the above coupling.
An application of the monotone convergence theorem implies that
\begin{align}\label{eqn:2236}
\cW\big(\cC\big(\mvy^{[n]}, t; \ell\big)\big)\convas\cW\big(\cC\big(\mvy,t; \ell\big)\big)\, ,
\ \ \text{ as }\ \ n\to\infty\, .
\end{align}
Using \eqref{eqn:2234} in conjunction with \eqref{eqn:2235} and \eqref{eqn:2236} and passing to the $n\to\infty$ limit, it follows that 
$\big(f_{\ell}\big(G(\mvy, t)\big), t\geq 0\big)$ is a submartingale.
}

To prove (\ref{eq:submart_toprove}), note that
\begin{align}\label{eqn:3}
\cA f_{\ell}(G)
&
=\sum_{i\neq k}z_kz_i\Bigg(\frac{(z_k+z_i)^2}{\sigma_2+2z_kz_i}
-\frac{z_k^2}{\sigma_2}\Bigg)\\
&\hskip25pt
+\sum_{\substack{i\geq 1\\ i\neq k }}\sum_{\substack{j>i\\ j\neq k}}z_iz_j\Bigg(\frac{z_k^2}{\sigma_2+2z_iz_j}
-\frac{z_k^2}{\sigma_2}\Bigg)
=:T_1+T_2.\notag
\end{align}
Now
\begin{align*}
T_1
&
=\sum_{i\neq k}\frac{z_kz_i}{\sigma_2\big(\sigma_2+2z_kz_i\big)}
\bigg(\big(z_k^2+z_i^2+2z_kz_i\big)\sigma_2-z_k^2\big(\sigma_2+2z_kz_i\big)\bigg)\\
&\hskip25pt
=\sum_{i\neq k}\frac{z_kz_i}{\sigma_2\big(\sigma_2+2z_kz_i\big)}
\bigg(z_i^2\sigma_2+2z_kz_i\big(\sigma_2-z_k^2\big)\bigg)\\
&\hskip50pt
=z_k\sum_{i\neq k}\frac{z_i^3}{\sigma_2+2z_kz_i}
+2z_k^2\sum_{i\neq k}\frac{z_i^2\big(\sum_{j\neq k} z_j^2\big)}{\sigma_2\big(\sigma_2+2z_kz_i\big)}.
\end{align*}
Using $2z_kz_i\leq\sigma_2$, we get
\begin{align}\label{eqn:4}
T_1
\geq
z_k\sum_{i\neq k}\frac{z_i^3}{\sigma_2+2z_kz_i}
+\frac{z_k^2}{\sigma_2^2}\times\bigg(\sum_{i\neq k} z_i^2\bigg)^2.
\end{align}
Next
\begin{align}\label{eqn:5}
-T_2
=
z_k^2\sum_{\substack{i\geq 1\\ i\neq k}}\sum_{\substack{j>i\\ j\neq k}}\frac{2z_i^2z_j^2}{\sigma_2\big(\sigma_2+2z_iz_j\big)}
\leq
z_k^2\sum_{\substack{i\geq 1\\ i\neq k}}\sum_{\substack{j>i\\ j\neq k}}
\frac{2z_i^2z_j^2}{\sigma_2^2}
\leq\frac{z_k^2}{\sigma_2^2}\times\bigg(\sum_{i\neq k} z_i^2\bigg)^2.
\end{align}
Combining \eqref{eqn:3}, \eqref{eqn:4}, and \eqref{eqn:5} yields \eqref{eq:submart_toprove}.

Now take $\ell=1$ and assume that
$f_1(G(\mvy,0)) = y_1^2/\sigma_2(\mvy)>1-\eps$.
Define the stopping time
\[T=\inf\big\{s\geq 0: f_1(G(\mvy,s))\leq 4/5\big\},\]
where \ch{the} infimum of an empty set is understood to be $+\infty$.
Since  $f_1(G(\mvy,t))$ is a bounded submartingale,  $f_1(G(\mvy,\infty)):=\lim_{t\to\infty}f_1(G(\mvy,t))$ exists almost surely.
Thus,
\[1-\eps
\leq\E\big[f_1(G(\mvy,T))\big]
\leq\frac{4}{5}\bP\big(f_1(G(\mvy,T))\leq 4/5\big)+\bP\big(f_1(G(\mvy,T))>4/5\big).\]
This shows that $\bP(T=\infty)\geq 1-5\eps$.
Then \eqref{eqn:leader-doesnot-change} follows upon noting that $\cW(\cC(\mvy,t; 1))>2\max_{i\geq 2}\cW(\cC_i(\mvy,t))$ if $f_1(G(\mvy,t))>4/5$.
\qed

\medskip

\noindent\ch{{\bf Proof of Proposition \ref{prop:4.1}:}
Combined with Proposition~\ref{prop:submartingale}, the bound given in Proposition \ref{prop:initial-condition-case-i} below completes the proof.
\qed}
\begin{prop}\label{prop:initial-condition-case-i}
Under Condition I, for all $\eps>0$, there exists $\lambda>0$ depending only on $\eps$ and the sequence $\{\mvx^{\sss(n)}\}_{n\geq 1}$ such that for all $n \ge 1$,
\begin{align*}
\bP\Bigg(\frac{\big(x_1^{\sss(n)}(t_{\lambda})\big)^2}{\sigma_2\big(\mvx^{\sss(n)}(t_{\lambda})\big)}\geq 1-\eps\Bigg)\geq 1-\eps\, ,
\end{align*}
where $t_{\lambda}=t_\lambda(\mvx^{(n)})=\lambda+1/\sigma_2(\mvx^{\sss(n)})$.
\end{prop}
\noindent{\bf Proof:}
First consider the particular sequences $\mvy^{(n)}$ defined in \eqref{eq:errg-weights}, for which $\sigma_2(\mvy^{(n)})=n^{-1/3}$. We then have $t_\lambda(\mvy^{(n)})=\lambda + n^{1/3}$.
Ignoring vertex weights and multiple edges, the graph $G(\mvy^{(n)},\lambda+n^{1/3})$ has the same law as the Erd\H{o}s-R\'enyi random graph on $n$ vertices, where edges are placed independently between each pair of vertices with probability
\[
p_{\lambda} =p_\lambda(n)= \bP\big(\mathrm{Poisson}\big(n^{-4/3}(\lambda+n^{1/3})\big)\geq 1\big) = (1+o(1))\left(\frac{1}{n}+\frac{\lambda}{n^{4/3}}\right)\, .
\]
It thus follows from \cite[Theorem A.1]{janson08susceptibility} that
\begin{equation}\label{eq:special_case}
\liminf_{n \to \infty}\ \bP\Bigg(\frac{\big(y_1^{\sss(n)}(\lambda+n^{1/3})\big)^2}{\sigma_2\big(\mvy^{\sss(n)}(\lambda+n^{1/3})\big)}\geq 1-\eps\Bigg)\geq 1-\eps\, ,
\end{equation}
provided $\lambda$ is sufficiently large.

Fix any sequence $\{\mvx^{\sss(n)}\}_{n\geq 1}$ satisfying Condition I.
Using \eqref{eqn:weight-tau-con}, we see that with $t_{\lambda}=t_{\lambda}(\mvx^{\sss(n)})$,
\begin{gather}
\mvx^{\sss(n)}(t_\lambda) \convd \mvzeta(\lambda) \text{ on } \ldown,\ \text{ and} \label{eqn:weight-converge}\\
\mvy^{\sss(n)}(\lambda + n^{1/3}) \convd \mvzeta(\lambda) \text{ on } \ldown.\label{eqn:34}
\end{gather}
Now fix $\lambda$ sufficiently large so that \eqref{eq:special_case} holds.
Then using \eqref{eqn:weight-converge} and \eqref{eqn:34},
\begin{align*}
\liminf_{n \to \infty}\ \bP\Bigg(\frac{\big(x_1^{\sss(n)}(t_{\lambda})\big)^2}{\sigma_2\big(\mvx^{\sss(n)}(t_{\lambda})\big)}> 1-2\eps\Bigg)& \geq
\bP\left(\frac{|\gamma_1(\lambda)|^2}{\sigma_2(\mvzeta(\lambda))} > 1-2\eps \right) \\
& \ge \bP\left(\frac{|\gamma_1(\lambda)|^2}{\sigma_2(\mvzeta(\lambda))} \ge 1-\eps \right) \\
& \ge
\limsup_{n \to \infty}\ \bP\Bigg(\frac{\big(y_1^{\sss(n)}(\lambda+n^{1/3})\big)^2}{\sigma_2\big(\mvy^{\sss(n)}(\lambda+n^{1/3})\big)}\geq 1-\eps\Bigg) \geq  1-\eps\, ,
\end{align*}
as desired.
\qed

\medskip

It remains to prove tightness of the family of random variables in Proposition \ref{prop:4.2}.
We will need some properties of the distributional limits $(\mvzeta(\lambda), \lambda \in \bR)$.
By \cite[Proposition 18 and Equation (80)]{aldous-limic}, for any $\eps>0$, we can choose $\lambda_{\eps}>0$ such that with probability at least $1-\eps$, the following three assertions hold simultaneously:
\begin{align}\label{eqn:1818}
\lambda_{\eps}\sigma_2\big(\mvzeta(-\lambda_{\eps})\big)\geq 1/2,\ \ \
\Bigg|\frac{1}{\sigma_2\big(\mvzeta(-\lambda_{\eps})\big)}-\lambda_{\eps}\Bigg|\leq 1,\ \text{ and }\
\lambda_{\eps}\times|\gamma_1(-\lambda_{\eps})|\leq\eps.
\end{align}
Writing $t_{\lambda}=\lambda+1/\sigma_2(\mvx^{\sss(n)})$ as before, and using \eqref{eqn:weight-converge},  we can
choose $n_0(\eps)$ such that for all $n\geq n_0(\eps)$,
\begin{align}\label{eqn:18}
\lambda_{\eps}\sigma_2\big(\mvx^{\sss(n)}\big(t_{-\lambda_{\eps}}\big)\big)\geq 1/2,\ \ \
\Bigg|\frac{1}{\sigma_2\big(\mvx^{\sss(n)}\big(t_{-\lambda_{\eps}}\big)\big)}-\lambda_{\eps}\Bigg|\leq 1,\ \text{ and }\
\lambda_{\eps}\times x_1^{\sss(n)}\big(t_{-\lambda_{\eps}}\big)\leq\eps
\end{align}
with probability at least $1-2\eps$.

We will also use a simple lemma leveraging the coupling between $(\mvx(t),t \ge 0)$ and $(G(\mvx,t),t \ge 0)$. Recall that $\cC(\mvx,t;j)$ is the component containing vertex $j$ in $G(\mvx,t)$.
\begin{lem}\label{lem:expected-cluster-size}
For $\mvx=(x_i,i \ge 1) \in \ldown$ and $t>0$, if $t\sigma_2(\mvx)<1$ then for all $j \ge 1$,
\[\E\big[\cW\big(\cC(\mvx,t;j)\big)\big]\leq\frac{x_j}{1-t\sigma_2(\mvx)}.\]
\end{lem}
\noindent{\bf Proof:}
By the definition of $G(\mvx,t)$, for each $i \ne j$, the probability that there is at least one edge between $i$ and $j$ in $G(\mvx,t)$ is $1-\exp(-tx_ix_j)$, and these events are mutually independent. It follows that for $i\neq j$,
\begin{align*}
\bP\big(i\in\cC(\mvx,t;j)\big)
&\leq\sum_{k\geq 0}\sum_{j_1\neq j}\ldots\sum_{j_k\neq j}(tx_j x_{j_1})\cdot(tx_{j_1}x_{j_2})\cdot\ldots\cdot(tx_{j_{k-1}}x_{j_k})\cdot(tx_{j_k}x_{i})\\
&=tx_ix_j\sum_{k\geq 0}\big(t\big(\sigma_2(\mvx)-x_j^2\big)\big)^k
=\frac{tx_ix_j}{1-t\big(\sigma_2(\mvx)-x_j^2\big)}\, .
\end{align*}
Hence
\begin{align}\label{eqn:3333}
\E\big[\cW(\cC(\mvx,t;j))\big]
&=x_j+\sum_{i\neq j}x_i\cdot\bP\big(i\in\cC(\mvx,t;j)\big)
\leq x_j+\frac{tx_j \big(\sigma_2(\mvx)-x_j^2\big)}{1-t\big(\sigma_2(\mvx)-x_j^2\big)}
\notag\\
&=\frac{x_j}{1-t\big(\sigma_2(\mvx)-x_j^2\big)}
\leq 
\frac{x_j}{1-t \sigma_2(\mvx)}\, ,
\end{align}
as desired.
\qed

\medskip

\noindent{\ch{\bf Proof of Proposition \ref{prop:4.2}:}}
Let us first describe the core idea in words.
Using the above estimates and lemma, we will show that the maximal component from time $t_{-\lambda_{\eps}}$ is unlikely to become very large by time $t_{-2}$.
On the other hand, by \eqref{eqn:weight-converge}, the maximal component at time $t_{-2}$ is reasonably large.
This implies that with high probability the leader changes at least once in the time interval $[t_{-\lambda_\eps} , t_{-2}]$. We now make this idea precise.

Let $\mvx^{\sss(n)}\big(t_{-\lambda_{\eps}}\big)=:\mvz=(z_i,i \ge 1)$ and suppose we start the multiplicative coalescent with $\mvz$ as the initial configuration. 
Run this process for $\lambda_\eps -2$ units of time.
Let $E(n)$ be the event given in \eqref{eqn:18}. 
Using Lemma \ref{lem:expected-cluster-size}, on the event $E(n)$,
\begin{align*}
&\E\big[\cW\big(\cC(\mvz,\lambda_{\eps}-2; 1)\big)\ \big|\ \mvz\big]
\leq\frac{z_1}{1-(\lambda_{\eps}-2)\sigma_2(\mvz)}\\
&\hskip60pt
=
\frac{z_1}{\sigma_2(\mvz)}\times\frac{1}{\frac{1}{\sigma_2(\mvz)}-\lambda_{\eps}+2}
\leq\frac{z_1}{\sigma_2(\mvz)}\leq\frac{\eps}{\lambda_{\eps}\sigma_2(\mvz)}
\leq 2\eps.
\end{align*}
Hence, for all $\delta>0$ and $n\geq n_0(\eps)$, by an application of Markov's inequality,
\begin{align}\label{eqn:19}
\bP\big(\cW\big(\cC\big(\mvz,\lambda_{\eps}-2; 1\big)\big)\geq\delta\big)\leq \bP(E(n)^c)+2\eps/\delta \le 2\eps\big(1+1/\delta\big).
\end{align}
However by \eqref{eqn:weight-converge} for all $\eta>0$, there exists $\delta_{\eta}\in (0,1)$ and $n_1(\eta)\geq 1$ such that for all $n\geq n_1(\eta)$,
\begin{align}\label{eqn:20}
\bP\big(x_1^{\sss(n)}\big(t_{-2}\big)\geq\delta_{\eta}\big)\geq 1-\eta.
\end{align}
Taking $\eps=\eta\delta_{\eta}$ and
combining \eqref{eqn:19} and \eqref{eqn:20}, we see that for all $n\geq n_0(\eps)\vee n_1(\eta)$,
\[\bP\left(\bigg(L(\mvx^{(n)}) - \frac{1}{\sigma_2(\mvx^{(n)})}\bigg)^{\mathlarger{-}}
\geq
\lambda_{\eps}\right)\leq\eta+2\eps(1+1/\delta_{\eta})\le 5\eta\, ,\]
which completes the proof. 
\qed

\section{Proof of Theorem \ref{thm:number-of-changes}}
The proof of Proposition \ref{prop:4.2} in fact shows that for any $\eta>0$, there is $\lambda_{\eta}^{\sss(1)}>0$ such that for all $n$ sufficiently large,
\[
\bP\big(\text{The leader changes between times }t=t_{-\lambda_{\eta}^{\sss(1)}}\text{ and }t=0\big)\geq 1-\eta.
\]
By repeating the same argument, we can choose
$\lambda_{\eta}^{\sss(1)}<\lambda_{\eta}^{\sss(2)}<\lambda_{\eta}^{\sss(3)}\ldots$
such that for all $j\geq 1$, for all $n$ sufficiently large,
\[
\bP\Big(\text{The leader changes between times }t=t_{-\lambda_{\eta}^{\sss(j+1)}}\text{ and }t=t_{-\lambda_{\eta}^{\sss(j)}}\Big)\geq 1-\eta/2^{j}.
\]
It thus follows that for any $j\geq 1$ and $\eta>0$,
\[
\limsup_n\ \bP\big(N^{\sss(n)}\leq j\big)\leq 2\eta,
\]
which completes the proof.

\section{Proof of Theorem~\ref{thm:main_ii}}
\label{sec:main_ii_proof}
In this section, we work under Condition II.
Thus, throughout this section, $\mvc = (c_i,i \ge 1)$ satisfies $c_i \asymp i^{-\alpha}$ for some fixed $\alpha \in (1/3,1/2)$, and $\big(\mvx^{(n)},\, n \ge 1\big)$ satisfies \eqref{eq:main_conv} and \eqref{eqn:555} with this limiting sequence $\mvc$.

Define $\phi:[0,\infty)\to\bR$ and $\Phi_{\lambda}:[0,\infty)\to\bR$ for each $\lambda>0$ by setting
\[\phi(s)=\sum_{i\geq 1}c_i^2\Big(\frac{c_i s-1+e^{-c_i s}}{c_i s}\Big)\ \ \text{ and }\ \
 \Phi_{\lambda}(s):=\lambda s-s\phi(s).\]
\ch{Since $u\mapsto (u-1+e^{-u})$ is convex, it follows that $s\phi(s)$ is also convex, and consequenctly $\Phi_{\lambda}$ is concave.}
Also observe that for each $\lambda>0$, $\Phi_{\lambda}$ has a unique positive zero,
which we denote by $s_0(\lambda)$.
Further, $s_0(\lambda)$ is strictly increasing on $(0,\infty)$, and $s_0(\lambda)\uparrow\infty$ as $\lambda\uparrow\infty$.

Recall the definitions of the processes $W_\lambda(\cdot)$ and $\overline{W}_\lambda(\cdot)$ from \eqref{eqn:wvs-def} and \eqref{eqn:bar-lam-def}, and of the vector $\mvzeta(\lambda):=(|\gamma_{i}(\lambda)|,\ i\geq 1)\in\ldown$ of ranked excursion lengths of $\overline{W}_\lambda$. Note that for all $\lambda > 0$,
\begin{equation}\label{eq:wphi-diff}
W_\lambda(s) - \Phi_\lambda(s) = \sum_{i \ge 1} c_i\big(\ind_{\{\xi_i \le s\}} - \bP(\xi_i \le s)\big)\, .
\end{equation}
The next theorem bounds the length and squared sum of jump sizes of the longest excursion $\gamma_1(\lambda)$.
\begin{thm}\label{lem:largest-component-big}
For all $\eta > 0$, there exist positive constants $\lambda_0$ and $C$ depending only on $\eta$ and on $(c_i, i\geq 1)$ such that
for all $\lambda\geq\lambda_0$,
\begin{align}\label{eqn:17}
\bP\Big(\big|\gamma_{1}(\lambda)\big|\leq
(1-\eta) s_0(\lambda)
\Big)
\leq
\exp\Big(-C\lambda^{\frac{1}{2(1-2\alpha)}}\log\log\lambda\Big)~,
\end{align}
and
\begin{equation}\label{eq:size_lbd}
\bP\Big(
\big|\gamma_{1}(\lambda)\big|\ge
(1+\eta) s_0(\lambda)
\Big)
\le
C\lambda^{-\frac{1}{1-2\alpha}}~.
\end{equation}
Further, there exist positive constants $\delta$, $C'$ and $\lambda'$ depending only on $\mvc$ such that for all $\lambda\geq\lambda'$,
\begin{align}\label{eqn:16}
\bP\Big(\sum_{i\in\gamma_{1}(\lambda)}c_i^2\leq (1+2\delta)\lambda\Big)\leq
C'\lambda^{-\frac{1}{1-2\alpha}}.
\end{align}
\end{thm}
Theorem~\ref{lem:largest-component-big} will imply rather straightforwardly that the bound of Proposition~\ref{prop:initial-condition-case-i} also holds under Condition II; this is the content of the next proposition.
\begin{prop}\label{prop:initial-condition-case-ii}
Under Condition II, for all $\eps>0$, there exists $\lambda_0>0$ depending only on $\eps$ and the sequence $\{\mvx^{\sss(n)}\}_{n\geq 1}$ such that for all $\lambda\geq\lambda_0$,
\begin{align*}
\liminf_{n\to\infty}\ \bP\Bigg(\frac{\big(x_1^{\sss(n)}(t_{\lambda})\big)^2}{\sigma_2\big(\mvx^{\sss(n)}(t_{\lambda})\big)}\geq 1-\eps\Bigg)\geq 1-\eps\, ,
\end{align*}
where $t_{\lambda}=t_\lambda(\mvx^{(n)})=\lambda+1/\sigma_2(\mvx^{\sss(n)})$.
\end{prop}
Combined with Proposition~\ref{prop:submartingale}, Theorem~\ref{thm:main_ii} follows immediately. Thus, the rest of the section is devoted to the proofs of Theorem~\ref{lem:largest-component-big} and
Proposition~\ref{prop:initial-condition-case-ii}.
%
%

We first list some simple asymptotics which will be useful in the proof.
Note that $u-1+e^{-u}\asymp u$ on $u\geq 1$, and 
$u-1+e^{-u}\asymp u^2$ on $u\in [0, 1]$.
For $s>0$, let 
\begin{align}\label{eqn:def-i-0}
	i_0(s)=\min\big\{i\geq 1\ :\ c_i s< 1\big\}.
\end{align}
\ch{Note that $i_0(s)\geq 2$ when $s\geq 1/c_1$, and}
\begin{align}\label{eqn:46}
s\phi(s)
=
\sum_{i=1}^{i_0(s)-1} c_i\big(c_i s-1+e^{-c_i s}\big)+\sum_{i=i_0(s)}^{\infty} c_i\big(c_i s-1+e^{-c_i s}\big)
\asymp
\sum_{i=1}^{i_0(s)-1} c_i(c_i s)+\sum_{i=i_0(s)}^{\infty} c_i(c_i s)^2
\end{align}
on $s\geq 1/c_1$.
From the relation $c_i\asymp i^{-\alpha}$ it follows that 
\begin{align}\label{eqn:777}
i_0(s)\asymp s^{1/\alpha}\ \ \text{ on }\ \  s\geq 1/c_1.
\end{align}
Combining \eqref{eqn:46} with \eqref{eqn:777} and the relation $c_i\asymp i^{-\alpha}$, we get
\begin{align}\label{item:3}
s\phi(s) = \sum_{i=1}^{\infty} c_i\big(c_i s-1+e^{-c_i s}\big)\asymp s^{(1-\alpha)/\alpha}\ \ \text{ on }\ \ s\geq 1/c_1.
\end{align}
Similarly, using the relations 
\begin{align}\label{eqn:666}
1-e^{-u}\asymp 1\ \text{ on }\ u\geq 1\  \text{ and }\ 1-e^{-u}\asymp u\  \text{ on }\ u\in[0,1],
\end{align} 
we get
\begin{align}\label{item:2}
\sum_{i=1}^{\infty} c_i^2\big(1-e^{-c_i s}\big)\asymp s^{(1-2\alpha)/\alpha}\ \ \text{ on }\ \ s\geq 1/c_1~,
\end{align}
and, for any $k\geq 2$ and $j\geq 1$,
\begin{equation}
\sum_{i=j}^{\infty} c_i^k\big(1-e^{-c_i s}\big)\asymp s j^{1-\alpha(k+1)}\ \text{ on }\ 0\leq s\leq 1/c_j ~.
\label{item:1}
\end{equation}
Finally, note that $1-e^{-u}-ue^{-u}\geq 0$ for $u\in\bR$, and further,  
$1-e^{-u}-ue^{-u}\asymp 1$ on $u\geq 1$ and 
$1-e^{-u}-ue^{-u}\asymp u^2$ on $u\in[0,1]$.
Hence, 
\begin{equation}\label{eqn:44}
\sum_{i=1}^{\infty} c_i\big(1-e^{-c_i s}-c_i s e^{-c_i s}\big)\asymp s^{(1-\alpha)/\alpha}\ \ \text{ on }\ \ s\geq 1/c_1.
\end{equation}
For fixed $\eta\in (0,1)$ and $w\in[1-\eta, 1]$,
\begin{align}\label{eqn:45}
\big|\sum_{i=1}^{\infty} c_i\big(c_i s e^{-c_i s w}-c_i s e^{-c_i s }\big)\big|
\leq \sum_{i=1}^{\infty} c_i^2 s\big(1-e^{-c_is\eta}\big)\ \ \ch{\text{ for }\ \  s\geq 0}\, .
\end{align}
If $c_1 s\eta\geq 1$, then
\begin{align}\label{eqn:47}
\sum_{i=1}^{\infty} c_i^2 s\big(1-e^{-c_is\eta}\big)
=
\sum_{i=1}^{i_0(s\eta)-1} c_i^2 s\big(1-e^{-c_is\eta}\big)
+
\sum_{i=i_0(s\eta)}^{\infty} c_i^2 s\big(1-e^{-c_is\eta}\big)
\leq C s^{(1-\alpha)/\alpha}\eta^{(1-2\alpha)/\alpha}~,
\end{align}
where $i_0(\cdot)$ is as in \eqref{eqn:def-i-0}, and the last inequality uses \eqref{eqn:666}, \eqref{eqn:777}, and the \ch{  relation} $c_i\asymp i^{-\alpha}$.
If $c_1 s\eta< 1$, then 
\begin{align}\label{eqn:48}
\sum_{i=1}^{\infty} c_i^2 s\big(1-e^{-c_is\eta}\big)\leq s^2\eta\sum_{i=1}^{\infty} c_i^3
\leq s^{(1-\alpha)/\alpha}\cdot (c_1\eta)^{(1-3\alpha)/\alpha}\cdot\eta\cdot\sum_{i=1}^{\infty} c_i^3
\leq C s^{(1-\alpha)/\alpha}\eta^{(1-2\alpha)/\alpha}.
\end{align}
Combining \eqref{eqn:45}, \eqref{eqn:47}, and \eqref{eqn:48} with \eqref{eqn:44}, we see that for $\eta > 0$ sufficiently small, uniformly over $w \in [1-\eta,1]$,
\begin{equation}
\sum_{i=1}^{\infty} c_i\big(1-e^{-c_i s}-c_i s e^{-c_i sw}\big)\asymp s^{(1-\alpha)/\alpha}\ \ \text{ on }\ \ s\geq 1/c_1.
\label{item:4}
\end{equation}

For $u > 0$, define
\begin{align}\label{eqn:def-Z-u-1}
Z_u^{\sss(1)}=\sup_{s\le u}\Big|\sum_{i\geq 1}c_i\big(\ind_{\{\xi_i\leq s\}}-\bP(\xi_i\leq s)\big)\Big|
= \sup_{s \le u} |W_\lambda(s) - \Phi_\lambda(s)|\, ,
\end{align}
where the second equality follows from \eqref{eq:wphi-diff}.
Note that
\begin{align*}
Z_u^{\sss(1)} & =
\lim_{k \to \infty}
\sup_{s \le u: ks \in \bN}
\Big|\sum_{i\geq 1}c_i\big(\ind_{\{\xi_i\leq s\}}-\bP(\xi_i\leq s)\big)\Big|
\\
& =
\lim_{k \to \infty}\ \Big(\lim_{n \to \infty}
\sup_{s \le u~:~ ks \in \bN}~
\Big|\sum_{i=1}^n c_i\big(\ind_{\{\xi_i\leq s\}}-\bP(\xi_i\leq s)\big)\Big|\Big)\, ,
\end{align*}
where the second step uses the fact that $\sum_{i\ge 1}c_i\big(\ind_{\{\xi_i\leq s\}} -\bP(\xi_i\leq s)\big)$ is almost surely convergent for any fixed $s$ \ch{(which is a simple consequence of Kolmogorov's
three-series theorem \cite{durrett2019probability} and the fact that $\mvc\in l^3$)}.
Defining
\[
Z_{u,n}^{\sss(1)} =
\sup_{s\le u}\ \Big|\sum_{i= 1}^n c_i\big(\ind_{\{\xi_i\leq s\}}-\bP(\xi_i\leq s)\big)\Big|\, ,
\]
for all $k$ we have
$\sup_{s \le u:\ ks \in \bN}
\Big|\sum_{i=1}^n c_i\big(\ind_{\{\xi_i\leq s\}}-\bP(\xi_i\leq s)\big)\Big|\le Z_{u,n}^{\sss(1)}$.
Combined with the above identities and an application of Fatou's lemma, this yields that for all $z$,
\begin{equation}\label{eq:zreduce}
\bP\bigg(Z_u^{\sss(1)}> z\bigg) \le
\sup_{n \ge 1}\ \bP\bigg(Z_{u,n}^{\sss(1)}> z\bigg)\, \ch{.}
\end{equation}
The point of this is that it will shortly allow the use of discrete-time martingale bounds which {\em a priori} only hold for observables of finite-dimensional random variables.
We likewise define
\begin{align}\label{eqn:def-Z-u-2}
Z_u^{\sss(2)} =\sup_{s\le u}\Big|\sum_{i\geq 1}c_i^2\big(\ind_{\{\xi_i\leq s\}}-\bP(\xi_i\leq s)\big)\Big|\ \ \text{ and }\
\
Z_{u,n}^{\sss(2)} =
\sup_{s\le u}\Big|\sum_{i= 1}^n c_i^2\big(\ind_{\{\xi_i\leq s\}}-\bP(\xi_i\leq s)\big)\Big|\, .
\end{align}
The same argument as above shows that for all $z$,
\begin{equation}\label{eq:zreduce2}
\bP\Big(Z_u^{\sss(2)}> z\Big) \le
\sup_{n \ge 1}\ \bP\Big(Z_{u,n}^{\sss(2)}> z\Big)\, .
\end{equation}
\begin{lem}\label{lem:concentration}
There exists a constant $C_{\ref{lem:concentration}}$ depending only on $(c_i, i\geq 1)$ such that for all $u \ge 1/C_{\ref{lem:concentration}}$ and all $x \ge 1/C_{\ref{lem:concentration}}$,
\begin{gather}
\bP\Big(Z_u^{\sss(1)}\geq x u^{(1-2\alpha)/(2\alpha)}\Big)
\leq \exp\Big(-C_{\ref{lem:concentration}} x\log\log x\Big),\ \
\text{ and}\label{eqn:14}\\
\bP\Big(Z_u^{\sss(2)}\geq x \Big)
\leq \exp\Big(-C_{\ref{lem:concentration}} x\log\log x\Big).\label{eqn:15}
\end{gather}
\end{lem}
\noindent{\bf Proof:}
By \eqref{eq:zreduce}, it suffices to prove \eqref{eqn:14} with $Z_{u,n}^{\sss(1)}$ in place of $Z_u^{\sss(1)}$. We consider positive and negative fluctuations separately, defining
\begin{align*}
U_{u,n}= \sup_{s\le u}\Big(\sum_{i= 1}^n c_i\big(\ind_{\{\xi_i\leq s\}}-\bP(\xi_i\leq s)\big)\Big),\ \text{ and }\
V_{u,n}= \sup_{s\le u}\Big(\sum_{i= 1}^n c_i\big(\bP(\xi_i\leq s)-\ind_{\{\xi_i\leq s\}}\big)\Big).
\end{align*}
Note that $Z_{u,n}^{\sss(1)}\leq U_{u,n}+V_{u,n}$.

We shall apply Theorem~\ref{thm:empirical_tail} to the random variable $u^{-(1-2\alpha)/(2\alpha)} U_{u,n}$.
To this end, observe that for $s\le u$,
\begin{align*}
 \var\bigg(u^{-(1-2\alpha)/(2\alpha)}
\sum_{i=1}^n c_i\big(\ind_{\ch{\{\xi_i \le s\}}} - \bP(\xi_i \le s)\big)\bigg)
& = u^{-(1-2\alpha)/\alpha}\sum_{i=1}^n c_i^2 \bP(\xi_i \le s)\big(1-\bP(\xi_i \le s)\big)
\\
& \leq u^{-(1-2\alpha)/\alpha}\sum_{i\ge 1} c_i^2(1-e^{-c_i u}) \le K'\, ,
\end{align*}
where the last bound is due to (\ref{item:2}).
Now Lemma~\ref{lem:bound-expectation} given below implies that 
$\E U_{u,n} \le K_{\ref{lem:bound-expectation}}u^{(1-2\alpha)/(2\alpha)}$ provided $u \ge 1/c_1$. 
Using this, an application of Theorem~\ref{thm:empirical_tail} with $V=K'$ and $A=c_1 u^{-(1-2\alpha)/(2\alpha)}$ gives
\begin{align*}
& \bP\bigg( \frac{U_{u,n}}{u^{(1-2\alpha)/(2\alpha)}} \ge K_{\ref{lem:bound-expectation}}+x
\bigg)\le
\exp\left(
 - \frac{xu^{(1-2\alpha)/(2\alpha)}}{4c_1} \log\bigg(1 + 2\log\Big(1+\frac{xu^{(1-2\alpha)/(2\alpha)}}{c_1K'}\Big)\bigg)
\right)\, .
\end{align*}
\ch{(Note that the supremum in the definition of $U_{u,n}$ can be equivalently taken over the countable set $[0, u]\cap\bQ$, where $\bQ$ denotes the set of rationals. 
Thus, Theorem \ref{thm:empirical_tail} is applicable here.)}
An identical treatment yields a similar tail bound for $V_{u,n}$, which, combined with the preceding tail bound, establishes \eqref{eqn:14}.

The proof of \eqref{eqn:15} is very similar to the proof of \eqref{eqn:14}. The key difference comes from the improved bound
\[\ch{
\var\big(\sum_{i=1}^n c_i^2\big(\ind_{\{\xi_i \le s\}} - \bP(\xi_i \le s)\big)\big) = 
\sum_{i=1}^n \var\big(c_i^2\big(\ind_{\{\xi_i \le s\}} - \bP(\xi_i \le s)\big)\big) 
\leq \sum_{i=1}^{\infty}c_i^4<\infty\, ,}
\]
and the use of \eqref{eqn:889} rather than \eqref{eqn:888}. This explains the gain of a factor $u^{(1-2\alpha)/(2\alpha)}$ in the bound of \eqref{eqn:15} relative to that of \eqref{eqn:14}. We omit the details to avoid repetition.
\qed

\begin{lem}\label{lem:bound-expectation}
There exists a constant $K_{\ref{lem:bound-expectation}}$ depending only on $(c_i, i\geq 1)$ such that for all $u \ge 1/c_1$ and all $n \ge 1$,
\begin{gather}
\E\bigg(
\sup_{s\le u}\ch{\Big|\sum_{i= 1}^n c_i\big(\ind_{\{\xi_i\leq s\}}-\bP(\xi_i\leq s)\big)\Big|}
\bigg)\leq K_{\ref{lem:bound-expectation}} u^{(1-2\alpha)/(2\alpha)},\mbox{ and }\label{eqn:888}\\
\E
\bigg(
\sup_{s\le u}\ch{\Big|\sum_{i= 1}^n c_i^2\big(\ind_{\{\xi_i\leq s\}}-\bP(\xi_i\leq s)\big)\Big|}
\bigg)\leq K_{\ref{lem:bound-expectation}}\, .\label{eqn:889}
\end{gather}
\end{lem}
\noindent{\bf Proof:}
Let $(\xi_1',\ldots,\xi_n')$ be an independent copy of $(\xi_1,\ldots,\xi_n)$.
Writing $\E_{\mvxi}[\cdot]:=\E[\cdot | \xi_1,\ldots,\xi_n]$, we have\ch{
\begin{align*}
& \E\bigg(
\sup_{s\le u}\Big|\sum_{i= 1}^n c_i\big(\ind_{\{\xi_i\leq s\}}-\bP(\xi_i\leq s)\big)\Big|
\bigg)\\
&\hskip20pt=
\E\Big(\sup_{s\le u}\Big| \E_{\mvxi}\Big(\sum_{i= 1}^n c_i\big(\ind_{\{\xi_i\leq s\}}
-\ind_{\{\xi_i'\leq s\}}\big)\Big)\Big|\Big)\\
&\hskip40pt
\leq
\E\E_{\mvxi}\Big(\sup_{s\le u}\Big|\sum_{i= 1}^n c_i\big(\ind_{\{\xi_i\leq s\}}
-\ind_{\{\xi_i'\leq s\}}\big)\Big|\Big)
=
\E\Big(\sup_{s\le u}\Big|\sum_{i= 1}^n c_i\big(\ind_{\{\xi_i\leq s\}}
-\ind_{\{\xi_i'\leq s\}}\big)\Big|\Big)\, .
\end{align*}}
Introducing i.i.d. random variables $\eps_1,\ldots,\eps_n$ with $\bP(\eps_1=1)=1/2=\bP(\eps_1=-1)$,
we see that\ch{
\begin{align}\label{eqn:27}
& \E\bigg(
\sup_{s\le u}\Big|\sum_{i= 1}^n c_i\big(\ind_{\{\xi_i\leq s\}}-\bP(\xi_i\leq s)\big)\Big|
\bigg)\\
&\hskip20pt
\leq
\E\Big(\sup_{s\le u}\Big|\sum_{i= 1}^n c_i\eps_i\big(\ind_{\{\xi_i\leq s\}}
-\ind_{\{\xi_i'\leq s\}}\big)\Big|\Big)\notag\\
&\hskip40pt
\leq
\E\Big(\sup_{s\le u}\Big|\sum_{i= 1}^n c_i\eps_i\ind_{\{\xi_i\leq s\}}\Big|\Big)
+\E\Big(\sup_{s\le u}\Big|\sum_{i= 1}^n c_i\big(-\eps_i\big)\ind_{\{\xi_i'\leq s\}}\Big|\Big)\notag\\
&\hskip60pt
=
2\E\Big(\sup_{s\le u}\Big|\sum_{i= 1}^n c_i\eps_i\ind_{\{\xi_i\leq s\}}\Big|\Big)
=
2\E\E_{\mvxi}\Big(\sup_{s\le u}\Big|\sum_{i= 1}^n c_i\eps_i\ind_{\{\xi_i\leq s\}}\Big|\Big).\nonumber
\end{align}}
To bound the inner expectation, let $\{k(1),\ldots,k(M)\} \subseteq \{1,\ldots,n\}$ be the indices $k$ for which $\xi_k \le u$, listed in increasing order.
Then\ch{
\[
\sup_{s \le u} \Big|\sum_{i= 1}^n c_i\eps_i\ind_{\{\xi_i\leq s\}}\Big| =
\max_{1\leq j\leq M} \Big|\sum_{i=1}^j c_{k(i)} \eps_{k(i)}\Big|.
\]}
Write  $S(j) = \sum_{i=1}^j c_{k(i)} \eps_{k(i)}$.
Then, conditionally on $\mvxi$, $(S(j),~ 0 \le j \le M)$ is a martingale.
Hence,
\begin{align*}
&\ch{\E_{\mvxi}\Big(\sup_{s\le u}\Big|\sum_{i= 1}^n c_i\eps_i\ind_{\{\xi_i\leq s\}}\Big|\Big)
=}
\E_{\mvxi}\Big[\max_{1\leq j\leq M}|S(j)|\Big]
\leq
\Big(\E_{\mvxi}\Big[\max_{1\leq j\leq M}S(j)^2\Big]\Big)^{1/2}\\
&\hskip50pt
\leq
2\Big(\E_{\mvxi}\big[S(M)^2\big]\Big)^{1/2}
=
2\Big(\E_{\mvxi}\Big[\big(\sum_{i=1}^n c_i\eps_i\ind_{\ch{\{\xi_i\leq u\}}}\big)^2\Big]\Big)^{1/2}
\leq
2\Big(\sum_{i\geq 1} c_i^2\ind_{\ch{\{\xi_i\leq u\}}}\Big)^{1/2},
\end{align*}
where the second step uses the Cauchy-Schwarz inequality, the third step uses Doob's inequality, and
\ch{the final inequality follows by directly expanding the square and using the independence of $\eps_i$, $1\leq i\leq n$}.
Combined with \eqref{eqn:27}, this gives
\begin{align*}
&\ch{\E\bigg(
\sup_{s\le u}\Big|\sum_{i= 1}^n c_i\big(\ind_{\{\xi_i\leq s\}}-\bP(\xi_i\leq s)\big)\Big|
\bigg)}
\le
4\E\Big[\Big(\sum_{i\geq 1} c_i^2\ind_{\{\xi_i\leq u\}}\Big)^{1/2}\Big]\notag
\\
&\hskip50pt
\le
4\Big(\E\big[\sum_{i\geq 1} c_i^2\ind_{\{\xi_i\leq u\}}\big]\Big)^{1/2}
=
4\Big(\sum_{i\geq 1} c_i^2\big(1-e^{-c_i u}\big)\Big)^{1/2}
\le K_{\ref{lem:bound-expectation}} u^{(1-2\alpha)/(2\alpha)}\, ,
\end{align*}
the final inequality holding by (\ref{item:2}). 
This yields \eqref{eqn:888}.
The proof of \eqref{eqn:889} is essentially identical and is omitted.
\qed

\medskip

The final step before the proofs of Theorem~\ref{lem:largest-component-big} and Proposition~\ref{prop:initial-condition-case-ii} is to collect a few facts about the functions $\phi$ and $s_0$ defined at the start of the section, and about the function $\psi:[0,\infty) \to \bR$ defined by
\[
\psi(s)=\sum_{i\geq 1}c_i^2\big(1-e^{-c_i s}\big).
\]
%
%
%
\begin{lemma}\label{lem:deterministic-functions}
The following properties are
satisfied by $\phi, \psi$, and $s_0$.
\begin{enumerate}[label=\upshape({\alph*})]
\item
$s_0(\lambda)\asymp \lambda^{\alpha/(1-2\alpha)}$ on $[s_0^{-1}(1/c_1),\ \infty)$.
\item
There exists $\eta_0\in(0, 1/2)$ such that for all $\eta\in(0,\eta_0]$,
\[\sup_{s\geq 2/c_1}\frac{\phi(\eta s)}{\phi(s)}
\leq \sup_{s\geq 2/c_1}\frac{\phi\big((1-\eta)s\big)}{\phi(s)}
=:1-\eps_{\eta}~,\]
where $\eps_{\eta}>0$ and $\lim_{\eta\to 0}\eps_{\eta}=0$.
\item
There exist $\delta_0>0$ and $\lambda_0>0$ such that for all $\lambda\geq\lambda_0$,
\[\psi\big(s_0(\lambda)\big)\geq (1+\delta_0)\lambda.\]
\item
For all $s>0$ and $\eta\in(0,1)$, $\psi\big(\eta s\big)\geq \eta\psi(s)$. Further, if
\[g(\eta):=\sup_{s\geq 1/c_1}\frac{\psi(\eta s)}{\psi(s)}~,\]
then $g(\eta)\to 0$ as $\eta\to 0$.
\end{enumerate}
\end{lemma}
\noindent{\bf Proof:}
Note that $\lambda=\phi(s_0(\lambda))\asymp \big(s_0(\lambda)\big)^{(1-2\alpha)/\alpha}$ on $\{\lambda : s_0(\lambda)\geq 1/c_1\}$, where the final step uses \eqref{item:3}.
This proves (a).

Next, a direct calculation shows that the function $u\mapsto (u-1+e^{-u})/u$ is increasing on $\bR$. This implies that
\[\sup_{s\in[2/c_1,\ \infty)}\frac{\phi(\eta s)}{\phi(s)}
\leq \sup_{s\in[2/c_1,\ \infty)}\frac{\phi\big((1-\eta)s\big)}{\phi(s)}\]
for $0<\eta<1/2$. Now
\begin{align}\label{eqn:36}
\frac{\phi\big(s\big)}{\phi((1-\eta)s)}-1
& = \frac{\sum_{i\geq 1}c_i\Big(\eta\big(1-e^{-c_i s}\big)+e^{-c_i s}-e^{-c_i(1-\eta)s}\Big)}{(1-\eta)s\cdot\phi((1-\eta)s)}
\notag\\
&\geq \frac{\sum_{i\geq 1}c_i\Big(\eta\big(1-e^{-c_i s}\big)+e^{-c_i s}-e^{-c_i(1-\eta)s}\Big)}{Ks^{(1-\alpha)/\alpha}}\notag\\ &
=
\frac{\sum_{i\geq 1}c_i\int_{1-\eta}^1\Big(1-e^{-c_i s}-c_i s e^{-c_i su}\Big)du}{Ks^{(1-\alpha)/\alpha}},
\end{align}
where the second step makes use of \eqref{item:3} and is valid whenever $c_1(1-\eta)s\geq 1$.
By \eqref{item:4}, there is $\eta_0\in(0, 1/2)$ such that
uniformly over $u \in [1-\eta_0,1]$,
\[
\sum_{i \ge 1} c_i (1-e^{-c_i s} - c_i s e^{-c_i s u}) \asymp s^{(1-\alpha)/\alpha}\ \text{ on }\ s\geq 1/c_1.
\]
This last observation combined with (\ref{eqn:36}) shows that \ch{there exists $K'>0$ such that for every $\eta \in (0,\eta_0]$,
\[
\frac{\phi\big(s\big)}{\phi((1-\eta)s)}-1 \geq K'\eta\ \text{ for all }\ s \ge 1/(c_1(1-\eta_0)).
\]
which completes the proof of (b).}

Since $1-e^{-u}\geq (u-1+e^{-u})/u$ for $u\in (0,\infty)$, it immediately follows that $\psi(s)\geq\phi(s)$ for $s\geq 0$.
Now
\[
s\big(\psi(s)-\phi(s)\big) =
\sum_{i \ge 1} c_i (1-e^{-c_i s} - c_i s e^{-c_i s})
 \asymp s^{(1-\alpha)/\alpha} \asymp s\phi(s)\ \text{ on }\  s\geq 1/c_1,
\]
where we use \eqref{item:4} and \eqref{item:3} for the final two asymptotic equivalences.
It follows that there exist $\delta_0>$ and $s_
\star>0$ such that
$\psi(s)\geq (1+\delta_0)\phi(s)$
for all $s\geq s_\star$, which in turn implies the existence of a $\lambda_0$ such that
\[\psi\big(s_0(\lambda)\big)\geq (1+\delta_0)\phi\big(s_0(\lambda)\big)=(1+\delta_0)\lambda\]
for all $\lambda\geq\lambda_0$. This proves (c).

Finally, since $1-\exp(-\eta s)\geq \eta(1-e^{-s})$, it follows that $\psi\big(\eta s\big)\geq \eta\psi(s)$.
Now note that
\begin{align*}
g(\eta)& \leq
\sup_{\eta c_1s\geq 1}\frac{\psi(\eta s)}{\psi(s)}+
\sup_{1\leq c_1s\leq 1/\eta }\frac{\psi(\eta s)}{\psi(s)}\\
& \leq K\Big(\eta^{(1-2\alpha)/\alpha}+\sup_{1\leq c_1s\leq 1/\eta }\frac{\eta s}{s^{(1-2\alpha)/\alpha}}\Big)
\leq K\Big(\eta^{(1-2\alpha)/\alpha}+\frac{\eta^{(1-2\alpha)/\alpha}}{c_1^{(3\alpha-1)/\alpha}}\Big),
\end{align*}
where the second step uses \eqref{item:2} and \eqref{item:1}. This completes the proof of (d).
\qed

\medskip

\noindent{\bf Proof of Theorem \ref{lem:largest-component-big}:}
Let $\delta_0$, $\lambda_0$, $\eta_0$, and $g(\cdot)$ be as in Lemma \ref{lem:deterministic-functions}.
Choose any $\eta\in(0,\eta_0]$ small enough that
\begin{align}\label{eqn:37}
\big(1-\eta-g(\eta)\big)(1+\delta_0)>1+\delta_0/2.
\end{align}
Now note that when $s_0(\lambda)\geq 2/c_1$,
\begin{align*}
\Phi_{\lambda}\big((1-\eta)s_0(\lambda)\big)
&=(1-\eta)s_0(\lambda)\Big(\lambda-\phi\big((1-\eta)s_0(\lambda)\big)\Big)\\
&\geq (1-\eta)s_0(\lambda)\Big(\lambda-(1-\eps_{\eta})\phi\big(s_0(\lambda)\big)\Big)
=(1-\eta)s_0(\lambda)\lambda\eps_{\eta}\geq K_{\eta}\lambda^{\frac{1-\alpha}{1-2\alpha}}
\end{align*}
for some $K_{\eta}>0$, where the second step uses Lemma \ref{lem:deterministic-functions} (b) and the last step uses Lemma \ref{lem:deterministic-functions}~(a).
Similarly, $\ \Phi_{\lambda}\big(\eta s_0(\lambda)\big)\geq K_{\eta}\lambda^{\frac{1-\alpha}{1-2\alpha}}\ $ when $s_0(\lambda)\geq 2/c_1$.
Since $\Phi_{\lambda}(\cdot)$ is concave, it follows that
\begin{align}\label{eqn:38}
\Phi_{\lambda}\big(s\big)\geq K_{\eta}\lambda^{\frac{1-\alpha}{1-2\alpha}},
\ \ \text{ for }\ \ \eta s_0(\lambda)\leq s\leq (1-\eta)s_0(\lambda).
\end{align}
Recall the definitions of $W_{\lambda}(\cdot)$ and $Z_u^{\sss(1)}$ from \eqref{eqn:wvs-def} and \eqref{eqn:def-Z-u-1} respectively.
Then
\begin{align}\label{eqn:39}
W_{\lambda}(s)\geq \Phi_{\lambda}(s)-Z_{s_0(\lambda)}^{\sss(1)}\ \ \text{ for }\ \ 0\leq s\leq s_0(\lambda).
\end{align}
Using \eqref{eqn:38}, \eqref{eqn:39}, Lemma \ref{lem:deterministic-functions}(a), and applying \eqref{eqn:14} with $u=s_0(\lambda)$ and
$x=\theta \lambda^{\frac{1}{2(1-2\alpha)}}$
where $\theta>0$ is very small, we see that for all $\lambda$ sufficiently large,
\begin{align}\label{eqn:40}
\bP\Big(\inf\big\{W_{\lambda}(s)\ :\ \eta s_0(\lambda)\leq s\leq (1-\eta)s_0(\lambda)\big\}>0\Big)
\geq 1-\exp\big(-Cf_1(\lambda)\big),
\end{align}
where $f_1(\lambda)=\lambda^{\frac{1}{2(1-2\alpha)}}\log\log\lambda$.
Write $\gamma_\star$ for $\gamma(\lambda; s_0(\lambda)/2)$-the excursion of the reflected process $\overline W_{\lambda}(\cdot)$ alive at time $s_0(\lambda)/2$.
Then \eqref{eqn:40} shows that $\gamma_{\star}$ is alive when $\eta s_0(\lambda)\leq s\leq (1-\eta)s_0(\lambda)$
with probability at least $1-\exp\big(-Cf_1(\lambda)\big)$.
Thus, with probability at least $1-\exp\big(-Cf_1(\lambda)\big)$, $|\gamma_1(\lambda)|\ge (1-2\eta)s_0(\lambda)$.
This proves \eqref{eqn:17}.

\ch{Let us make a note here of the following bound which we will use shortly: 
For all $\lambda$ sufficiently large,
\begin{equation}\label{eq:sizeupperbd}
\bP\Big(
\exists s \in \big[(1-\eta)s_0(\lambda)\, ,\ (1+\eta)s_0(\lambda)\big]\text{ such that }
W_\lambda(s) < 0
\Big)
\geq 1-\exp\big(-Cf_1(\lambda)\big)\, .
\end{equation}
To see this, observe that Lemma \ref{lem:deterministic-functions}~(b) implies that there exists $\eps'_{\eta}>0$ such that
$\phi\big(\big(1+\eta/2\big)s\big)\geq \big(1+\eps'_{\eta}\big)\phi(s)$ for all $s\geq 2/c_1$.
Consequently, for $\lambda\geq s_0^{-1}(2/c_1)$,
\begin{align*}
\Phi_{\lambda}\big(\big(1+\eta/2\big)s_0(\lambda)\big)
&=\big(1+\eta/2\big)s_0(\lambda)\Big(\lambda-\phi\big(\big(1+\eta/2\big)s_0(\lambda)\big)\Big)\\
&\leq \big(1+\eta/2\big)s_0(\lambda)\Big(\lambda-(1+\eps_{\eta}')\phi\big(s_0(\lambda)\big)\Big)
=-\big(1+\eta/2\big)s_0(\lambda)\lambda\eps_{\eta}'
\leq -K_{\eta}'\lambda^{\frac{1-\alpha}{1-2\alpha}}\, ,
\end{align*}
for some $K_{\eta}'>0$, where the last step uses Lemma \ref{lem:deterministic-functions}~(a).
Since $\Phi_{\lambda}$ is concave with $\Phi_{\lambda}(0)=0=\Phi_{\lambda}\big(s_0(\lambda)\big)$, $\Phi_{\lambda}$ is non-increasing on $[s_0(\lambda), \infty)$.
Hence, for $\lambda\geq s_0^{-1}(2/c_1)$,
$\Phi_{\lambda}(s)\leq -K_{\eta}'\lambda^{\frac{1-\alpha}{1-2\alpha}}$
for $\big(1+\eta/2\big)s_0(\lambda)\leq s\leq (1+\eta)s_0(\lambda)$.
Since
$W_{\lambda}(s)\leq \Phi_{\lambda}(s)+Z_{(1+\eta)s_0(\lambda)}^{\sss(1)}$ for $0\leq s\leq (1+\eta)s_0(\lambda)$, we can use \eqref{eqn:14} as before to get \eqref{eq:sizeupperbd}.
}

Now, for all sufficiently large $\lambda$,
\[\sum_{i\in\gamma_\star}c_i^2\geq Y:=\sum_{i\geq 1}c_i^2\ind\big\{\eta s_0(\lambda)\leq \xi_i\leq (1-\eta)s_0(\lambda)\big\}\]
with probability at least $1-\exp\big(-Cf_1(\lambda)\big)$.
Recall the definition of $Z_u^{\sss (2)}$ from \eqref{eqn:def-Z-u-2}.
Then
\[
Y\geq \psi\big((1-\eta)s_0(\lambda)\big)-\psi\big(\eta s_0(\lambda)\big)- 2Z_{s_0(\lambda)}^{\sss(2)}.
\]
It follows that for all sufficiently large $\lambda$, with probability at least $1-\exp\big(-Cf_1(\lambda)\big)$,
\begin{align}\label{eqn:41}
\sum_{i\in\gamma_\star}c_i^2
&\geq 
(1-\eta)\psi\big(s_0(\lambda)\big)-g(\eta)\psi\big(s_0(\lambda)\big)-2Z_{s_0(\lambda)}^{\sss(2)}\\
&\hskip50pt\geq 
\big(1-\eta-g(\eta)\big)(1+\delta_0)\lambda-2Z_{s_0(\lambda)}^{\sss(2)}
\geq \big(1+\delta_0/2\big)\lambda-2Z_{s_0(\lambda)}^{\sss(2)},\notag
\end{align}
where the first, second, and third inequality use Lemma \ref{lem:deterministic-functions}(d), Lemma \ref{lem:deterministic-functions}(c), and \eqref{eqn:37} respectively.
Using the tail bound \eqref{eqn:15} with $x=\theta\lambda$ (where $\theta>0$ is very small), we see that for all $\lambda$ sufficiently large,
\begin{align}\label{eqn:42}
\bP\Big(\sum_{i\in\gamma_\star}c_i^2\leq (1+4\delta)\lambda\Big)\leq
\exp\bigg(-K f_2(\lambda)\bigg)
\end{align}
where $\delta := \delta_0/16$ and  $f_2(\lambda):=\lambda\log\log\lambda$.

\ch{Let $\lambda_0>0$ be large enough so that all the preceding bounds hold for every $\lambda\geq\lambda_0$.
For the remainder of the proof, we work with a fixed $\lambda\geq\lambda_0$.
None of the constants in the subsequent argument will depend on $\lambda$.}

Recall from Section~\ref{sec:bfs-const} the breadth-first walk process $B^{\sss(n)}=(B_{\mvx^{\sss(n)},t_\lambda}(s),s \ge 0)$ associated with $G(\mvx^{\sss(n)},t_\lambda)$, and write
$\fC_\star^{\sss(n)}$ for $\fC^{\sss(n)}(\lambda; s_0(\lambda)/2)$-the component of $G(\mvx^{\sss(n)}, t_{\lambda})$ being explored by $B^{\sss(n)}$ at time $s_0(\lambda)/2$.
Also, let $\cF^{\sss (n)}$ and $\cR^{\sss (n)}$ respectively denote the $\sigma$-field generated by the process $B^{\sss(n)}$ up to the time when the exploration of $\fC_\star^{\sss(n)}$ concludes and the set of vertices found up to that time.

Let $E^{\sss(n)}$ be the event that the following happen:
\begin{inparaenumn}
\item
exploration of $\fC_\star^{\sss(n)}$ begins before time $\eta s_0(\lambda)$;
\item
exploration of $\fC_\star^{\sss(n)}$ ends between times
$(1-\eta) s_0(\lambda)$ and $(1+\eta) s_0(\lambda)$; and 
\item 
\[
\sum_{v\in\fC_\star^{\sss(n)}}\frac{(x_v^{\sss(n)})^2}{\sigma_2(\mvx^{\sss(n)})^2}
\geq\lambda\big(1+2\delta\big).
\]
\end{inparaenumn}
So in particular, if $E^{\sss(n)}$ occurs then
\begin{align}\label{eqn:55}
(1-2\eta)s_0(\lambda)\leq \cW(\fC_\star^{\sss(n)})\leq (1+\eta)s_0(\lambda),\ \text{ and }\
\cW\big(\cC(\mvx^{\sss(n)}, t_{\lambda}; v)\big)\leq 2\eta s_0(\lambda)\ \text{ for }\
v\in\cR^{\sss(n)}\setminus \fC_\star^{\sss(n)}.
\end{align}
Then \eqref{eqn:2929-I} and \eqref{eqn:2929-II} combined with \eqref{eqn:40}, \eqref{eq:sizeupperbd}, and \eqref{eqn:42}
imply that there exists $n_0=n_0(\lambda)$ such that 
\begin{align}\label{eqn:56}
\bP({E^{\sss(n)}})\ge 1-2\exp(-Cf_2(\lambda))\ \ \text{ for }\ \ n\geq n_0.
\end{align}

Let
\[
\mvx^{\star}=\mvx^{\star}(n) = \big(x_v^{\sss(n)},~ v \not \in \cR^{(n)}\big),
\]
and let
\[
\sigma_2^{\star}=\sigma_2^{\star}(n) =
\sigma_2(\mvx^{\sss (n)}) -
\sum_{v\in\cR^{(n)}}(x_v^{\sss(n)})^2.
\]
Defining $\lambda^{\star}$ by the identity $\lambda+1/\sigma_2(\mvx^{\sss (n)}) = -\lambda^{\star}+1/\sigma_2^{\star}$,
it follows that on the event $E^{\sss(n)}$,
\begin{align}\label{eq:enstar_bd}
\lambda^{\star}
=\frac{\sigma_2(\mvx^{\sss (n)})-\sigma_2^\star}{\sigma_2(\mvx^{\sss (n)})\sigma_2^\star}-\lambda
\geq\frac{\lambda(1+2\delta)(\sigma_2(\mvx^{\sss (n)}))^2}{\sigma_2(\mvx^{\sss (n)})\sigma_2^\star}-\lambda
\geq 2\delta\lambda\, .
\end{align}
Writing $\cS^{\sss(n)}$ for the set of components of $G(\mvx^{\sss n},t_\lambda)$ explored after exploring $\fC_{\star}^{\sss(n)}$,
we have
\begin{align}\label{eqn:99}
\sum_{\cC \in \cS^{\sss (n)}}  \cW(\cC)^2
=
 \sum_{v\notin\cR^{(n)}}x_v^{\sss (n)}\cdot  \cW\big(\cC(\mvx^{\sss (n)}, t_{\lambda}; v)\big).
\end{align}
Write $t_{-2\delta\lambda}^\star = -2\delta\lambda+ 1/\sigma_2^\star$.
On account of \eqref{eq:enstar_bd}, on the event $E^{\sss(n)}$,
$t^{\star}_{-2\delta\lambda}\geq t_{\lambda}$.
Consequently, for $v\notin\cR^{(n)}$,
\begin{align}\label{eqn:100}
&\ind_{E^{\sss(n)}}\cdot\E\Big[\cW\Big(\cC\big(\mvx^{\sss (n)}, t_{\lambda}; v\big)\Big)~\big|~\cF^{\sss (n)}\Big]
=
\ind_{E^{\sss(n)}}\cdot\E\Big[\cW\Big(\cC\big(\mvx^{\star}, t_{\lambda}; v\big)\Big)~\big|~\cF^{\sss (n)}\Big]
\\
&\hskip60pt\leq 
\E\Big[\cW\Big(\cC\big(\mvx^\star,t^\star_{-2\delta\lambda}; v\big)\Big)~\big|~\cF^{\sss (n)}\Big]
\leq 
\frac{x_v^{\sss(n)}}{1-t^\star_{-2\delta\lambda}\sigma_2^\star}
=
\frac{x_v^{\sss(n)}}{2\delta\lambda\sigma_2^\star}~,\notag
\end{align}
where the third step uses Lemma \ref{lem:expected-cluster-size}. 
Combining \eqref{eqn:99} and \eqref{eqn:100}, we get
\begin{align}\label{eqn:whoknows}
\E\Big[\sum_{\cC \in \cS^{\sss (n)}} \cW(\cC)^2\cdot\ind_{E^{\sss(n)}}~\big|~\cF^{\sss (n)}\Big]
\leq 
\sum_{v\notin\cR^{(n)}}x_v^{\sss (n)}\cdot \frac{x_v^{\sss (n)}}{2\delta\lambda\sigma_2^\star}
=\frac{1}{2\delta\lambda}~.
\end{align}
It thus follows that
\[\ind_{E^{\sss(n)}}\cdot\bP\Big(\sum_{\cC\in\cS^{\sss(n)}}\cW(\cC)^2\geq x
\ \Big|\ \cF^{\sss (n)}\Big)\leq \frac{1}{2\delta x\lambda}~.\]

Taking $x=\theta'\lambda^{\frac{2\alpha}{1-2\alpha}}$ (where $\theta'>0$ is very small) and using Lemma \ref{lem:deterministic-functions}~(a) together with \eqref{eqn:55} and \eqref{eqn:56}, we see that for all large $n$,
\begin{align}\label{eqn:57}
\bP\big(\fC_{\star}^{\sss(n)}\	\text{ is }\ \cC_1(\mvx^{\sss(n)},t_\lambda)\big)
\geq
1-\ch{2\exp\big(-Cf_2(\lambda)\big)-C'\lambda^{-\frac{1}{1-2\alpha}}}\, .
\end{align}
(Recall that $\cC_1(\mvx^{\sss(n)},t_\lambda)$ denotes the component of $G(\mvx^{\sss(n)},t_\lambda)$ having the largest mass.)
Using \eqref{eqn:57}, \eqref{eqn:55}, and \eqref{eqn:56}, it follows that for large $n$,
\begin{align*}
\bP\big( \cW\big(\cC_1(\mvx^{\sss(n)},t_\lambda)\big)> (1+\eta)s_0(\lambda) \big)
\le
\ch{4\exp(-Cf_2(\lambda)) + C'\lambda^{-\frac{1}{1-2\alpha}}}\, .
\end{align*}
Since \eqref{eqn:weight-tau-con} implies that $\cW(\cC_1(\mvx^{\sss(n)},t_\lambda)) \convd |\gamma_1(\lambda)|$,
this proves (\ref{eq:size_lbd}).

Finally, by \cite[Lemma 5.5]{bhamidi-hofstad-sen}, as $n\to\infty$,
\begin{align}\label{eqn:22}
\sum_{v\in\cC_1(\mvx^{\sss(n)},t_\lambda)}\frac{(x_v^{\sss(n)})^2}{\sigma_2(\mvx^{\sss(n)})^2}
\convd
\sum_{i\in\gamma_{1}(\lambda)}c_i^2.
\end{align}
Thus,
\begin{align*}
&\bP\Big(\sum_{i\in\gamma_{1}(\lambda)}c_i^2
<
 \lambda\big(1+2\delta\big)\Big)
\leq
\liminf_n\ \bP\Big(\sum_{v\in\cC_1(\mvx^{\sss(n)},t_\lambda)}\frac{(x_v^{\sss(n)})^2}{\sigma_2(\mvx^{\sss(n)})^2}< \lambda\big(1+2\delta\big)\Big)\\
&\hskip50pt
\leq \liminf_n\Big[ \bP\Big(\sum_{v\in\fC_{\star}^{\sss(n)}}\frac{(x_v^{\sss(n)})^2}{\sigma_2(\mvx^{\sss(n)})^2}\leq \lambda\big(1+2\delta \big)\Big)
+\bP\big(\fC_{\star}^{\sss(n)}\neq \cC_1(\mvx^{\sss(n)},t_\lambda)\big)\Big]
\leq C\lambda^{-\frac{1}{1-2\alpha}},
\end{align*}
where the last step follows from \eqref{eqn:56} and \eqref{eqn:57}.
This completes the proof of \eqref{eqn:16}.
\qed

\medskip

\noindent{\bf Proof of Proposition \ref{prop:initial-condition-case-ii}:}
Write $\sigma_2=\sigma_2(\mvx^{\sss(n)})$ for simplicity, and recall that $t_{\lambda}=\lambda+1/\sigma_2$.
Recall that $\cC_i(\mvx^{\sss(n)},t_\lambda)$ is the component of $\cG(\mvx^{\sss(n)}, t_{\lambda})$ with the $i$-th largest mass $x_i^{\sss(n)}(t_\lambda)$.

Fix $\eps>0$.
Using \eqref{eqn:weight-tau-con}, \eqref{eqn:22}, Theorem~\ref{lem:largest-component-big}, and the asymptotic from Lemma~\ref{lem:deterministic-functions}~(a), we can choose $\lambda(\eps)$ positive and large such that $\liminf_{n\to\infty}\pr\big(\cB_n(\lambda)\big)> 1-\eps$ for all $\lambda\geq\lambda(\eps)$, where
\begin{align}\label{eqn:21}
\cB_n(\lambda):=\Big\{x^{\sss(n)}_1\big(t_{\lambda}\big)\geq C_{\ref{eqn:21}}\lambda^{\frac{\alpha}{1-2\alpha}}\Big\}
\bigcap
\Big\{
\sum_{v\in\cC_1(\mvx^{\sss(n)},t_{\lambda})}\!\!\frac{(x^{\sss(n)}_v)^2}{\sigma_2^2} \geq (1+\delta)\lambda\Big\}.
\end{align}
From now on we work with a fixed $\lambda\geq\lambda(\eps)$ and $n$ large so that $\pr\big(\cB_n(\lambda)\big)\geq 1-\eps$.

We now reprise the argument leading to \eqref{eqn:whoknows}.
Define
\[\hat{\mvx}=\hat{\mvx}^{\sss (n)}=\big(x_v,\ v\notin\cC_1(\mvx^{\sss(n)},t_{\lambda})\big), \text{ and }\
\hat\sigma_2=\sigma_2(\hat{\mvx}^{\sss (n)})=\sigma_2-\!\!\sum_{v\in\cC_1(\mvx^{\sss(n)},t_{\lambda})}x_v^2.\]
We note here that $\cB_n(\lambda)$ is measurable w.r.t. the sigma-field generated by $\hat\mvx$.
Further, $x^{\sss(n)}_1(t_{\lambda})=\sigma_1(\mvx^{\sss(n)})-\sigma_1(\hat\mvx^{\sss(n)})$.

Define $\lambda'$ by the identity
$\lambda+1/\sigma_2=-\lambda'+1/\hat{\sigma}_2$.
Then on the event $\cB_n(\lambda)$,
\begin{align}\label{eqn:23}
\lambda'
=\frac{\sigma_2-\hat\sigma_2}{\sigma_2\hat\sigma_2}-\lambda
\geq\frac{(1+\delta)\sigma_2^2\lambda}{\sigma_2\hat\sigma_2}-\lambda
\geq\delta\lambda.
\end{align}
Write $\hat t_{-\delta\lambda}=-\delta\lambda+1/\hat\sigma_2$.
By \eqref{eqn:23}, on the event $\cB_n(\lambda)$, 
\begin{align}\label{eqn:61}
	\hat t_{-\delta\lambda}\geq t_{\lambda}.
\end{align}
Note that
\begin{align}\label{eqn:24}
\E\Big[\sigma_2\Big(\hat{\mvx}(\hat{t}_{-\delta\lambda})\Big)\ \big|\ \hat\mvx\Big]
&=\sum_{v\notin
\cC_1(\mvx^{\sss(n)},t_{\lambda})
}x_v^{\sss(n)}
\E\Big[\cW\Big(\cC\big(\hat\mvx,\hat t_{-\delta\lambda}; v\big)\Big)\ \big|\ \hat\mvx\Big]\notag\\
&\leq\sum_{v\notin
\cC_1(\mvx^{\sss(n)},t_{\lambda})
}
x_v^{\sss(n)}\frac{x_v^{\sss(n)}}{1-\hat t_{-\delta\lambda}\hat\sigma_2}
=\frac{1}{\delta\lambda},
\end{align}
where the second step uses Lemma \ref{lem:expected-cluster-size}.

For $z\geq 0$, define the event
$F(z) = \{\hat{x}_1(t_{\lambda}) \le z\}$, i.e., $F(z)$ denotes the event that all components of $\cG(\hat{\mvx},t_{\lambda})$ have weight at most $z$.
When $\cB_n(\lambda)$ occurs, by \eqref{eqn:61}, $\sigma_2\big(\hat{\mvx}(\hat{t}_{-\delta\lambda})\big) \ge (\hat{x}_1(t_{\lambda}))^2$.
Thus, on the event $\cB_n(\lambda)$, for any $z>0$,
\begin{align}\label{eqn:25}
1-\bP\big(F(z)\ \big|\ \hat\mvx\big)
\le \bP\Big(\sigma_2\big(\hat{\mvx}(\hat{t}_{-\delta\lambda})\big)>z^2\ \big|\ \hat\mvx\Big)
\le \frac{1}{z^2\delta\lambda}\ ,
\end{align}
where the final step uses \eqref{eqn:24}.

Now observe that, conditional on $\hat\mvx$, the graph $\cG(\mvx^{\sss (n)},t_{\lambda})\setminus\cC_1(\mvx^{\sss(n)},t_{\lambda})$ has the same distribution as $\cG(\hat{\mvx},t_{\lambda})$ conditional on $\hat\mvx$ and the event $F\big(x^{\sss(n)}_1(t_{\lambda})\big)$.
Hence, on the event $\cB_n(\lambda)$,
\begin{align*}
&\bP\Big(\sum_{i\geq 2}x_i^{\sss(n)}(t_{\lambda})^2\geq 1\ \big|\ \hat\mvx\Big)
=
\bP\Big(\sigma_2\big(\hat{\mvx}(t_{\lambda})\big)\geq 1\
\big|\ \hat\mvx,\ F\big(x^{\sss(n)}_1(t_{\lambda})\big)\Big)\\
&\hskip45pt\leq
\frac{\bP\Big(\sigma_2\big(\hat{\mvx}(t_{\lambda})\big)\geq 1 \big|\ \hat\mvx\Big)}{\bP\Big(F\big(x^{\sss(n)}_1(t_{\lambda})\big)\big|\ \hat\mvx\Big)}
\leq
\frac{\bP\Big(\sigma_2\big(\hat{\mvx}(\hat{t}_{-\delta \lambda})\big)\geq 1 \big|\ \hat\mvx\Big)}{\bP\Big(F\big(C_{\ref{eqn:21}}\lambda^{\frac{\alpha}{1-2\alpha}}\big)\big|\ \hat\mvx\Big)}
\leq \frac{1/(\delta\lambda)}{1-1/(\delta\lambda)},
\end{align*}
where the third step uses \eqref{eqn:61}, and
the last step uses \eqref{eqn:25} and presumes that $\lambda$ is large enough so that $C_{\ref{eqn:21}}\lambda^{\frac{\alpha}{1-2\alpha}}>1$.
Hence,
\begin{align*}
\bP\Big(\sum_{i\geq 2}x_i^{\sss(n)}(t_{\lambda})^2\geq 1\Big)
\leq \pr\big(\cB_n(\lambda)^c\big)+\frac{1/(\delta\lambda)}{1-1/(\delta\lambda)}
\leq \eps+\frac{1/(\delta\lambda)}{1-1/(\delta\lambda)}~.
\end{align*}
The rest is routine.
\qed

\section{Proof of Theorem \ref{thm:iv}}
\ch{The proof relies on the following lemma.
\begin{lem}\label{lem:234}
Suppose $\mvy=(y_i,\ i\geq 1)\in\ldown$ with $y_1>y_2>0$.
Let $y_1/y_2=(1+4\delta)$.
Fix $\eta>0$ such that
\begin{align}\label{eqn:2020}
\eta\leq 1/10\, ,\ \ \
\eta\sigma_2(\mvy)\leq \min\{1, y_1^2\}/4\, , \ \ \ \text{ and }\ \  \
1+\delta\leq \big(1+2\delta\big)\big(1-\eta\sigma_2(\mvy)\big)\, .
\end{align}
Suppose there exists $k\geq 1$ such that
\begin{align}\label{eqn:2021}
\bigg(\sum_{j=1}^k y_j\bigg)^2\geq \frac{\sigma_2(\mvy)}{\eta}\, \ \ \ \text{ and }\ \ \ 
\sum_{j\geq k+1}y_j^2\leq\eta\sigma_2(\mvy)\, .
\end{align}
Then 
\begin{align}\label{eqn:2022}
&\pr\big(\text{a change of leader does not occur in }G(\mvy, \cdot)\big)\notag\\
&\hskip40pt\geq
C_{\ref{eqn:2022}}\Big(1-\exp\big(-y_k^2\big)\Big)^k\times\exp\Big(-\big(1+2\delta\big)k^2 y_1^2\Big)
\times\left(\frac{\delta}{1+\delta}\right)^k=:g\big(\mvy, k,\delta\big)\, ,
\end{align}
where $C_{\ref{eqn:2022}}>0$ is an absolute constant.
\end{lem}
We will make use of the following result in the proof of Lemma \ref{lem:234}.
\begin{lem}[\cite{aldous-crit}, Lemma 20]\label{lem:235}
For $\mvz\in\ldown$ and $t_0\geq 0$, for all $s>\sigma_2(\mvz)$,
\[
\pr\big(\sigma_2\big(\mvz(t_0)\big)>s\big)\leq \frac{t_0 s\sigma_2(\mvz)}{s-\sigma_2(\mvz)}\, .
\]
\end{lem}
}

\noindent{\bf Proof of Lemma \ref{lem:234}:}\ch{
Let us first describe the proof idea in words.
Note that once the component of vertex $1$ achieves an `$l^2$-lead' over the other components in the sense of Proposition \ref{prop:submartingale}, then it stays the leader with positive probability for the rest of the process. 
So it is enough to show that with positive probability the component of vertex $1$ achieves this $l^2$-lead by some fixed time $t_0 > 0$, and it also stayed the leader in the time interval $[0, t_0]$.
We will now make this argument precise.
We can in fact take $t_0 = 1$, and will do so partway through the proof, but it is useful at the start to write $t_0$ to avoid confusion between time $1$ and the vertex label $1$.
}

\ch{
For $A\subseteq\{1, 2,\ldots\}$, $\ell\in A$, and $\mvu=\big(y_i; i\in A\big)$, write $\cC\big(\mvu, t_0; \ell\big)$ for the component of $\ell$ in the subgraph of $G(\mvy, t_0)$ that has $A$ as its vertex set.
(This notation agrees with the one introduced in Section \ref{sec:222}. 
The point we emphasize here is that all such components are coupled in that they are all subgraphs of $G(\mvy, t_0)$.)
Let
\begin{gather*}
\mvz_2=\big(y_i\, ;\ i=2\text{ or }i\geq k+1\big)\, ,\\
\mvz_j=\bigg(y_i\, ;\ i\in\{j\}\cup\big\{i\geq k+1\, :\, i\notin\cC\big(\mvz_{\ell}, t_0; \ell\big)
\text{ for }2\leq \ell\leq j-1 \big\}\bigg)\, ,\  \ j=3,\ldots,k\, ,\ \text{ and }\\
\mvz_{k+1}=\big(y_i\, ;\ i\geq k+1\text{ and }\ i\notin\cC\big(\mvz_{\ell}, t_0; \ell\big)
\text{ for }2\leq \ell\leq k\big)\, .
\end{gather*}
Define the events
\begin{gather*}
E_1:=\bigcap_{j=2}^k\big\{\text{there is an edge between }1\text{ and }j\text{ in }G(\mvy, t_0)\big\}\, ,\\
E_2:=\big\{\cW\big(\cC(\mvz_2, t_0; 2)\big)\leq (1+2\delta)y_2\big\}\, ,\\
E_j:=
\Bigg(\bigcap_{\ell=2}^{j-1}
\big\{\text{no edge between }j\text{ and any vertex in }\cC(\mvz_{\ell}, t_0; \ell)\big\}\Bigg)
\bigcap\big\{ \cW\big(\cC(\mvz_j, t_0; j)\big)\leq (1+2\delta)y_j \big\}
\end{gather*}
for $j=3,\ldots, k$, and 
\[
E_{k+1}:=\big\{\sigma_2\big(\mvz_{k+1}(t_0)\big)\leq 2\eta\sigma_2(\mvy)  \big\}\, .
\]
Then
\begin{align}\label{eqn:3334}
\pr(E_2)=1-\pr(E_2^c)\geq 1-\left(\frac{1}{(1+2\delta)y_2}\right)\left(\frac{y_2}{1-t_0\sum_{j\geq k+1}y_j^2}\right)
\geq \frac{\delta}{1+\delta}\, ,
\end{align}
where the second step uses Markov's inequality and the penultimate bound in \eqref{eqn:3333}, and the last step uses $t_0=1$, the second inequality in \eqref{eqn:2021}, and the third inequality in \eqref{eqn:2020}.
Similarly, for $j=3,\ldots, k$,
\begin{align}\label{eqn:3335}
\pr\big(E_j\, \big| E_2\cap\ldots\cap E_{j-1}\big)
\geq 
\exp\bigg(-y_j\sum_{\ell=2}^{j-1}\big(1+2\delta\big)y_{\ell}\bigg)\times
\left(\frac{\delta}{1+\delta}\right)
\geq 
\exp\bigg(-\big(1+2\delta\big)ky_1^2\bigg)\times
\left(\frac{\delta}{1+\delta}\right)\, .
\end{align}
Using Lemma \ref{lem:235} with $t_0=1$ and the second inequality in \eqref{eqn:2021},
\begin{align}\label{eqn:3336}
\pr\big(E_{k+1}\, |\, E_2\cap\ldots\cap E_k\big)
\geq 
1-\left(
\frac{2\eta\sigma_2(\mvy)\cdot\eta\sigma_2(\mvy)}{2\eta\sigma_2(\mvy)-\eta\sigma_2(\mvy)}
\right)
=1-2\eta\sigma_2(\mvy)\geq 1/2\, ,
\end{align}
where the last step uses the second inequality in \eqref{eqn:2020}.
Finally, $E_1$ is independent of the collection of events $\{E_2,\ldots, E_{k+1}\}$, and
\begin{align}\label{eqn:3337}
\pr(E_1)\geq \big(1-\exp(-y_k^2)\big)^k\, .
\end{align}
Combining \eqref{eqn:3334}, \eqref{eqn:3335}, \eqref{eqn:3336}, and \eqref{eqn:3337}, we get
\begin{align}\label{eqn:3338}
\pr\big(E_1\cap\ldots\cap E_{k+1}\big)\geq
\frac{1}{2}\Big(1-\exp\big(-y_k^2\big)\Big)^k\times\exp\Big(-\big(1+2\delta\big)k^2 y_1^2\Big)
\times\left(\frac{\delta}{1+\delta}\right)^k\, .
\end{align}
Now, on the event $E_1\cap\ldots\cap E_{k+1}$, the component of $1$ remains the unique leader all throughout the time interval $[0, t_0]$, and at time $t_0$, we have
\begin{align*}
\frac{\big(\cW\big(\cC(\mvy, t_0; 1)\big)\big)^2}{\sigma_2(\mvy(t_0))}
&\geq 
\frac{\big(\cW\big(\cC(\mvy, t_0; 1)\big)\big)^2}{\big(\cW\big(\cC(\mvy, t_0; 1)\big)\big)^2+2\eta\sigma_2(\mvy)}\\
&\geq\frac{\sigma_2(\mvy)/\eta}{\sigma_2(\mvy)/\eta+2\eta\sigma_2(\mvy)}
=\frac{1}{1+2\eta^2}\, ,
\end{align*}
where the first step uses the definitions of the events $E_1$ and $E_{k+1}$, and the second step uses the definition of the event $E_1$, and the first inequality in \eqref{eqn:2021}.
Using Proposition \ref{prop:submartingale}, 
\[
\pr\big(\text{the leader does not change in }\big(G(\mvy, t),\, t\geq t_0\big)\, \big|\, E_1\cap\ldots\cap E_{k+1} \big)
\geq\frac{1-8\eta^2}{1+2\eta^2}\geq\frac{1-8/100}{1+2/100}\, ,
\]
where the last step uses the fact $\eta\leq 1/10$.
Combining this with \eqref{eqn:3338} yields the claim. 
\qed
}

\vskip7pt

\noindent{\bf Completing the proof of Theorem \ref{thm:iv}:}
Consider $\mvx^{(n)}$, $n\geq 1$, satisfying Condition II with $c_1>c_2$.
Let $\theta\in(0, 1/3)$ be such that
\begin{align}\label{eqn:3339}
c_1(1-\theta)>c_2(1+\theta)\, .
\end{align}
Let $\mvzeta(\cdot)$ be the eternal multiplicative coalescent that satisfies \eqref{eqn:weight-tau-con} for each fixed $\lambda\in\bR$, and recall that
$\mvzeta(\lambda)=\big(|\gamma_1(\lambda)|, |\gamma_2(\lambda)|,\ldots\big)$.

Let $\eps=1/10$.
By \cite[Equation (65)]{aldous-limic}, we can choose $\lambda_{\eps}>0$ such that
\begin{align}\label{eqn:3340}
\pr\Big(
\forall\lambda\geq\lambda_{\eps},\  \
\lambda\cdot|\gamma_1(-\lambda)|\in \big(c_1(1-\theta),\, c_1(1+\theta)\big)\ \ \text{ and }\  \ 
\lambda\cdot|\gamma_2(-\lambda)|<c_2(1+\theta)
\Big)
>1-\eps\, .
\end{align}
Consider $\lambda_{\star}>\lambda_{\eps}$.
Let $\lambda_1>\lambda_2>\ldots>\lambda_\ell$ be such that
$\lambda_{\star}=\lambda_1$, $\lambda_{\eps}=\lambda_\ell$,
\begin{align}\label{eqn:3341}
\frac{c_1(1-\theta)}{\lambda_j}>\frac{c_2(1+\theta)}{\lambda_{j+1}}\, , \ \ \text{ and }\ \ 
\frac{2(1-\theta)}{\lambda_j}>\frac{(1+\theta)}{\lambda_{j+1}} \ \ \text{ for }\ \ j=1,\ldots, \ell-1\, .
\end{align}
(It is possible to choose $\lambda_1,\ldots, \lambda_\ell$ like this because of the relations \eqref{eqn:3339} and $2(1-\theta)>(1+\theta)$, the latter holding since $\theta<1/3$.)
By \eqref{eqn:weight-tau-con} and the Feller property of multiplicative coalescents \cite[Proposition 5]{aldous-crit},
$
\big(\mvx^{(n)}(t_{-\lambda_j})\, ,\ 1\leq j\leq \ell\big)\weakc
\big(\mvzeta(-\lambda_j)\, ,\ 1\leq j\leq \ell\big)
$,
which combined with \eqref{eqn:3340} implies that
\begin{align}\label{eqn:4040}
\liminf_{n \to \infty}\ \pr\big(\cA_n\big)\geq 1-\eps\, ,
\end{align}
where
$
\cA_n:=\bigcap_{j=1}^\ell\bigg
\{\lambda_j x_1^{(n)} \big(t_{-\lambda_j}\big)\in\big(c_1(1-\theta),\ c_1(1+\theta)\big) 
\ \text{ and }\ 
\lambda_j x_2^{(n)}\big(t_{-\lambda_j}\big)< c_2(1+\theta)
\bigg\}
$.
Assume that $\cA_n$ holds, and observe that if for some $\lambda\in[\lambda_2, \lambda_1)$, in the random graph $G\big(\mvx^{(n)}, t_{-\lambda}\big)$, the mass of a component other than the component containing the leader in $G\big(\mvx^{(n)}, t_{-\lambda_1}\big)$ becomes at least $c_1(1-\theta)/\lambda_1$,
then in $G\big(\mvx^{(n)}, t_{-\lambda_2}\big)$, either
{\bf (i)} that component will have merged with the component containg the leader in $G\big(\mvx^{(n)}, t_{-\lambda_1}\big)$, or
{\bf (ii)} that component and the component containing the leader in $G\big(\mvx^{(n)}, t_{-\lambda_1}\big)$ will remain disjoint.
However, \eqref{eqn:3341} shows that under {\bf (i)}, 
$x_1^{(n)}(t_{-\lambda_2})\geq 2c_1(1-\theta)/\lambda_1>c_1(1+\theta)/\lambda_2$,
and under {\bf (ii), }
$x_2^{(n)}(t_{-\lambda_2})\geq c_1(1-\theta)/\lambda_1>c_2(1+\theta)/\lambda_2$,
leading to a contradiction. 
Proceeding inductively, we see that on the event $\cA_n$, throughout the time interval $[t_{-\lambda_1},\, t_{-\lambda_\ell}]$, the component containing the leader in $G\big(\mvx^{(n)}, t_{-\lambda_1}\big)$ remains the unique leader.
Hence, \eqref{eqn:4040} implies that
\begin{align}\label{eqn:404}
\liminf_{n \to \infty}\ \pr\big(\cB_n\big)\geq 1-\eps\, ,
\end{align}
where $\cB_n$ is the event that there is a unique leader in $G\big(\mvx^{(n)}, t_{-\lambda_{\star}}\big)$, and the component containing the leader in $G\big(\mvx^{(n)}, t_{-\lambda_{\star}}\big)$ remains the unique leader throughout $[t_{-\lambda_{\star}},\, t_{-\lambda_{\eps}}]$.
(Note that the lower bound in \eqref{eqn:404} holds for all $\lambda_{\star}\in(\lambda_{\eps}, \infty)$.)

Since $\pr\big(|\gamma_2(-\lambda_{\eps})|>0\big)=1$, using \eqref{eqn:weight-tau-con} we see that 
\begin{align}\label{eqn:40404}
\liminf_{n \to \infty}\ \pr\big(x_2^{(n)}\big(t_{-\lambda_{\eps}}\big)>0\big)=1\, .
\end{align}
Let $\delta_n$ satisfy
$x_1^{(n)}\big(t_{-\lambda_{\eps}}\big)=(1+4\delta_n)x_2^{(n)}\big(t_{-\lambda_{\eps}}\big)$.
Since $\pr\big(|\gamma_1(-\lambda_{\eps})|>|\gamma_2(-\lambda_{\eps})|\big)=1$ 
\cite[Fact 1]{dhara-hofstad-leeuwaarden-sen-2},
\eqref{eqn:weight-tau-con} implies that
\begin{align}\label{eqn:3030}
\liminf_{n \to \infty}\ \pr\big(\delta_n>0\big)=1\, .
\end{align}
Let $\eta_n$ be the maximum of all $\eta$ that satisfy \eqref{eqn:2020} when $\mvy$ is replaced by $\mvx^{(n)}(t_{-\lambda_{\eps}})$ and $\delta$ is replaced by $\delta_n$.
Since $\sigma_2\big(\mvzeta(-\lambda_{\eps})\big)<\infty$ a.s. and $|\gamma_1(-\lambda_{\eps})|>0$ a.s., \eqref{eqn:weight-tau-con} and \eqref{eqn:3030} imply that 
\begin{align}\label{eqn:3031}
\liminf_{n \to \infty}\ \pr\big(\eta_n>0\big)=1\, .
\end{align}
Let $K_n$ be the minimum of all $k\geq 1$ that satisfy \eqref{eqn:2021} when $\mvy$ is replaced by $\mvx^{(n)}(t_{-\lambda_{\eps}})$ and $\eta$ is replaced by $\eta_n$ with the convention that $K_n=\infty$ if no such $k$ exists.
Since $\mvzeta(-\lambda_{\eps})\in\ldown\setminus l^1_{\downarrow}$ a.s., \eqref{eqn:weight-tau-con} and \eqref{eqn:3031} imply that we can choose $k_0\geq 1$ such that
\begin{align}\label{eqn:4041}
\liminf_{n \to \infty}\ \pr\big(K_n\leq k_0\big)> 1-\eps\, .
\end{align}
Using \eqref{eqn:weight-tau-con}, we see that
\begin{align}\label{eqn:4042}
g\big(\mvx^{(n)}(t_{-\lambda_{\eps}}), k_0, \delta_n\big)
\weakc
g\big(\mvzeta(-\lambda_{\eps}), k_0, \Delta\big)\, ,\ \ \text{ as }\ \ n\to\infty\, ,
\end{align}
where $g$ is as in \eqref{eqn:2022}, and $\Delta$ satisfies $|\gamma_1(-\lambda_{\eps})|=(1+4\Delta)\cdot |\gamma_2(-\lambda_{\eps})|$.
Since $\Delta>0$ a.s. and $|\gamma_{k_0}(-\lambda_{\eps})|>0$ a.s., the right side of \eqref{eqn:4042} is positive a.s.
Hence, we can choose $\kappa>0$ such that
\begin{align}\label{eqn:4043}
\liminf_{n \to \infty}\ 
\pr\big(g\big(\mvx^{(n)}(t_{-\lambda_{\eps}}), k_0, \delta_n\big)>\kappa\big)
> 1-\eps\, .
\end{align}

Define 
$
\cB_n':=\cB_n\cap\big\{ x_2^{(n)}\big(t_{-\lambda_{\eps}}\big)>0\big\}\cap\big\{K_n\leq k_0\big\}\cap\big\{ g\big(\mvx^{(n)}(t_{-\lambda_{\eps}}), k_0, \delta_n\big)>\kappa\big\}
$.
Combining \eqref{eqn:404}, \eqref{eqn:40404}, \eqref{eqn:4041}, and \eqref{eqn:4043}, we get
\begin{equation}\label{eqn:4044}
\liminf_{n \to \infty}\ \pr\big(\cB_n'\big)\geq 1-3\eps\, .
\end{equation}
Now
\begin{align}\label{eqn:444}
&\liminf_{n \to \infty}\ \pr\big(L(\mvx^{(n)})\leq t_{-\lambda_{\star}}\big)\notag\\
&\hskip15pt\geq
\liminf_{n \to \infty}\ 
\pr\Big(\cB_n'\cap\Big\{\text{a change of leader does not occur in }
\big(G(\mvx^{(n)}, t)\, ,\, t\geq t_{-\lambda_{\eps}}\big) \Big\}\Big)\notag\\
&\hskip30pt\geq
\liminf_{n \to \infty}\ 
\bE\big[\ind_{\cB_n'}\cdot g\big(\mvx^{(n)}(t_{-\lambda_{\eps}}), k_0, \delta_n\big)\big]
\geq
\kappa(1-3\eps)=7\kappa/10\, ,
\end{align}
where the second inequality uses Lemma \ref{lem:234} and the fact that $K_n\leq k_0$ on the event $\cB_n'$, and the third inequality uses \eqref{eqn:4044} and the fact that $g\big(\mvx^{(n)}(t_{-\lambda_{\eps}}), k_0, \delta_n\big)>\kappa$ on the event $\cB_n'$.
Since the right side of \eqref{eqn:444} is free of $\lambda_{\star}$, letting $\lambda_{\star}$ tend to infinity completes the proof.

\section*{Acknowledgements}
The authors thank James Martin and B\'{a}lazs R\'{a}th for helpful discussions about the results of  \cite{martin-rath}.
The authors also thank an anonymous referee for a careful reading of an earlier version of the paper and their detailed comments which led to significant improvements in the paper.
LAB was partially supported by NSERC Discovery Grant 341845.
SB was partially supported by NSF grants DMS-1613072, DMS-1606839 and ARO grant W911NF-17-1-0010.
SS was partially supported by a CRM-ISM fellowship, MATRICS grant MTR/2019/000745 from SERB, and by the Infosys Foundation, Bangalore.

\bibliographystyle{plainnat}
\bibliography{leader_bib}

\end{document}